\documentclass{amsart}
\usepackage[english]{babel}

\usepackage{amsmath}
\usepackage{amssymb}
\usepackage{amsfonts}
\usepackage{mathrsfs}
\usepackage{graphics}
\usepackage[matrix,arrow,curve]{xy}  
\usepackage{amscd} 

\usepackage{color}

\newtheorem{theorem}{Theorem}
\newtheorem{lemma}{Lemma}
\newtheorem{remark}{Remark}

\newtheorem{assum}{Assumption}

\newtheorem{defin}{Definition}

\newcommand{\ebox}{\fbox {} \smallskip}

\numberwithin{equation}{section}  
\numberwithin{theorem}{section}
\numberwithin{lemma}{section}
\numberwithin{remark}{section}
\numberwithin{prop}{section}
\numberwithin{assum}{section}
\numberwithin{coro}{section}
\numberwithin{defin}{section}

\def\{{\protect\lbrace}
\def\}{\protect\rbrace}

\def\const{{\rm{const}}}

\newcommand{\pa}{\partial}

\def\alf{\alpha}

\def\Del{\Delta}
\def\eps{\varepsilon}
\def\gam{\gamma}
\def\Gam{\Gamma}
\def\Gamo{\Gam_0}
\def\phi{\varphi}
\def\om{\omega}
\def\Om{\Omega}

\def\th{\theta}

\def\ov{\overline}
\def\ti{\tilde}

\def\wti{\widehat}

\newcommand{\bR}{\mathbb{R}}

\newcommand{\obR}{\ov \bR}

\newcommand{\V}{\vec V}
\newcommand{\W}{\vec W}
\newcommand{\F}{\Phi}
\newcommand{\FF}{\mathscr {F}}

\newcommand{\HH}{\mathsf {H}}

\newcommand{\gama}{\gam_{\alf}}
\newcommand{\gamap}{\gam_{\alf}^+}
\newcommand{\gamam}{\gam_{\alf}^-}
\newcommand{\gamapm}{\gam_{\alf}^{\pm}}
\newcommand{\gamamp}{\gam_{\alf}^{\mp}}

\newcommand{\ha}{h_{\alf}}

\newcommand{\yap}{y_{\alf}^+}
\newcommand{\yam}{y_{\alf}^-}
\newcommand{\yapm}{y_{\alf}^{\pm}}
\newcommand{\ye}{y_{\eps}}
\def\Gamo{\Gam_0}

\def\LA{L^{(a)}}
\def\LL{L^{(l)}}
\def\JA{J^{(a)}}
\def\JL{J^{(l)}}
\def\Del{\Delta}
\def\sh{\operatorname{sh}}
\def\ch{\operatorname{ch}}
\newcommand{\PR}{{\mathbb R}P}
\newcommand{\M}{Q}

\newcommand{\alfis}{\alf_{\textrm{is}}}
\renewcommand{\AA}{\mathscr {D}}

\begin{document}

\title[Geodesics in pseudo-Riemannian metrics]{On the local and global properties of geodesics in pseudo-Riemannian metrics}

\author{A.O. Remizov}
\address{CMAP, \'{E}cole Polytechnique CNRS, Route de Saclay, 91128 Palaiseau Cedex, France}
\email{alexey-remizov@yandex.ru}




\begin{abstract}
The paper is a study of geodesic in two-dimensional pseudo-Riemannian metrics.
Firstly, the local properties of geodesics in a neighborhood of generic parabolic points are investigated.
The equation of the geodesic flow has singularities at such points that leads to a curious phenomenon:
geodesics cannot pass through such a point in arbitrary tangential directions,
but only in certain directions said to be admissible (the number of admissible directions is generically 1 or 3).
Secondly, we study the global properties of geodesics in pseudo-Riemannian metrics
possessing differentiable groups of symmetries.
At the end of the paper, two special types of discontinuous metrics are considered.
\end{abstract}

\maketitle

\large

\section*{Introduction}

The paper presents a study of geodesics in two-dimensional pseudo-Riemannian metrics (the exact definitions are given below). A short section in the end of the paper is devoted to geodesics in metrics of two special types,
which are Riemannian at all points except for lines, where they are discontinuous.
The attention to such metrics is justified because they
have various applications in geometry, control theory and mechanics.

\medskip

By a pseudo-Riemannian metric on a manifold $\M$ we mean a quadratic form on the tangent bundle $T\M$,
whose signature may have different signs at different points of $\M$.
We shall always assume that the coefficients of pseudo-Riemannian metrics
smoothly ($C^{\infty}$) depend on a point of the manifold $\M$,
and change of the signature is only due to vanishing of the discriminant $\Del$.
The points of the set $\Del=0$ are often called {\it parabolic} points or {\it signature changing} points.
Generically, the set $\Del=0$ is a curve, which separates $\M$ into open domains with
constant signature: $(++)$ or $(+-)$ or $(--)$.
Pseudo-Riemannian metrics and their various geometric and physical aspects were considered by many authors; see e.g.
\cite{Hayward, K, KK2}. An example of a pseudo-Riemannian metric is the metric induced on a smooth surface embedded in three-dimensional Minkowski space \cite{GKhT, KhT, Rem-Pseudo}.

Let us note some specific properties of pseudo-Riemannian metrics.

For a non-parabolic point $q_0 \in Q$ and for any direction $p \in \PR$
there exists a unique geodesic passing through $q_0$ with given tangential direction $p$.
On the contrary, geodesics cannot pass through a parabolic point in arbitrary tangential directions,
but only in certain directions said to be admissible.
The number of admissible directions is generically 1 or 3, the number of geodesics with different admissible directions can also vary; see e.g. \cite{GR, Rem-Pseudo}.\footnote{
Similar results for three-dimensional pseudo-Riemannian metrics
were recently announced in \cite{PR11}.
}

The natural parametrization of geodesics in pseudo-Riemannian metrics
possess certain specific features.
By the natural parametrization we mean the parametrization of geodesics defined by the action functional
(this means that we treat them as extremal of the action functional).
Unlike the Riemannian case, the notion of natural parametrization does not completely coincide with the the arc-length parametrization, since the arc-length parametrization of isotropic geodesics does nor exist.

\medskip

The paper is organized as follows.
In Section~\ref{Definition} we discuss the definitions of geodesics through the action and the length functionals and establish the relationship between the natural and the arc-length parametrizations in pseudo-Riemannian metrics.

Section~\ref{Local} is devoted to the local properties of geodesics in a neighborhood of generic parabolic points
of pseudo-Riemannian metrics. This section consists of two parts. In the first part
we deal with unparametrized geodesics (in fact, we use a certain auxiliary parametrization, which is not the natural one). Most of the results of this section are taken from \cite{Rem-Pseudo}, but we reformulate them in a more convenient form. In the second part
we consider naturally parametrized geodesics and investigate their singularities at generic parabolic points.

In Section~\ref{Symmetries} we deal with metrics that possess differentiable groups of symmetries.
Namely, we consider metrics whose coefficients do not depend on one of the global coordinates on~$\M$,
and the group of symmetries consists of parallel translations of the coordinated plane.
Here a key point is the fact that the equation of unparametrized geodesics possesses a first integral,
which is often called the energy integral (the terminology comes from mechanics).
This allows us to reduce the second-order equation of unparametrized geodesics to a first-order implicit differential equation depending on a real parameter.
Using the qualitative theory of implicit differential equations~\cite{Dav-Japan},
we describe the global behavior of geodesics.
This constitutes the content of Section~\ref{Symmetries-Basic}.

Generalizations to metrics with another differentiable groups of symmetries are straightforward.
In general case, the existence of a first integral follows from Noether's theorem.
In the partial case, when the metric is induced on a surface of revolutions,
it also follows from well-known Clairaut's relation (see e.g.~\cite{Carmo}).
However, in many cases it is more convenient to operate with the the energy integral
using appropriate coordinates on $\M$ (as we have been doing in Section~\ref{Examples}).

In Section~\ref{R-and-L} we illustrate the obtained results for Riemannian and Lorenzian metrics.
By a Lorenzian metric we always mean a non-degenerate quadratic form on $T\M$
with the constant signature $(+-)$, but not necessarily constants coefficients (this terminology comes from physics).

Section~\ref{Pseudo-Riemannian} is the key part of the paper. Here we apply the obtained results
to pseudo-Riemannian metrics with differentiable groups of symmetries.
Although the obtained results seem rather simple, they allow to get some interesting consequences.
For instance, they allow to study in detail the global behavior of geodesics on surfaces of revolution
embedded in three-dimensional Euclidean or Minkowski space.
Section~\ref{Examples} presents two examples: geodesics on a sphere and geodesics on a torus in three-dimensional Minkowski space.

In the final Section~\ref{Klein-and-Grushin}, we investigate the global behavior of geodesics in
discontinuous metrics of two special types (called Klein type and Grushin type metrics)
when their coefficients do not depend on one of the coordinates.

\smallskip

\textbf{Acknowledgements}.
The research was supported by grant from FAPESP, proc. 2012/03960-2 for visiting ICMC-USP, S\~ao Carlos (Brazil).
I expresses deep gratitude to prof. Farid Tari for his attention to the work and useful discussions.
The last part of the paper was written during a visit to IHES in November--December 2012,
I am deeply grateful for the warm hospitality, excellent conditions and atmosphere.


\section{The definition of geodesics}\label{Definition}

Let $\M$ be a smooth two-dimensional manifold with the (local) coordinates $(x,y)$ and the
symmetrical covariant tensor field of the second order:
\begin{equation}
ds^2 = a(x,y) \, dx^2 + 2b(x,y) \, dx dy + c(x,y) \, dy^2
\label{1}
\end{equation}
said to be a metric on $\M$. We emphasize that \eqref{1} is not necessarily Riemannian (positive definite)
and even not necessarily non-degenerate at all points of $\M$.

For a curve $\gamma \colon I \to \M$ one can define the action and the length functionals:
\begin{equation*}
\JA(\gamma) = \int\limits_{\gamma} \bigl(a{\dot x}^2 + 2b{\dot x}\dot y + c{\dot y}^2\bigr)\,dt,
\quad
\JL(\gamma) = \int\limits_{\gamma} \sqrt{a{\dot x}^2 + 2b{\dot x}\dot y + c{\dot y}^2}\,dt,
\end{equation*}
where $\dot x = \frac{dx}{dt}$, $\dot y = \frac{dy}{dt}$.
Geodesics could be defined as extremals of the action functional, the
corresponding Euler-Lagrange equation reads
\begin{equation}
\frac{d}{dt} \LA_{\dot x} - \LA_{x} = 0, \ \ \ \frac{d}{dt}
\LA_{\dot y} - \LA_{y} = 0,
\label{2}
\end{equation}
where the Lagrangian $\LA(x,y,\dot x, \dot y) = a{\dot x}^2 +
2b{\dot x}\dot y + c{\dot y}^2$ is a function of the tangent bundle
$T\M$ homogeneous with respect to $(\dot x, \dot y)$. After straightforward
transformations, we get the system
\begin{equation}
\left \{ \
\begin{aligned}
& 2(a \ddot x + b \ddot y) = P, \quad
P = (c_x-2b_y) {\dot y}^2 - 2a_y {\dot x} {\dot y} -a_x {\dot x}^2, \\
& 2(b \ddot x + c \ddot y) = R, \quad
R = (a_y-2b_x) {\dot x}^2 - 2c_x {\dot x} {\dot y} - c_y {\dot y}^2, \\
\end{aligned}
\right.
\label{3}
\end{equation}
which can be written in the standard form
\begin{equation}
\left \{ \
\begin{aligned}
& \ddot x = \frac{cP-bR}{2\Del} = -\bigl(\Gamma^1_{11}{\dot x}^2 +
2\Gamma^1_{12}{\dot x}{\dot y} + \Gamma^1_{22}{\dot y}^2\bigr), \\
& \ddot y = \frac{aR-bP}{2\Del} = -\bigl(\Gamma^2_{11}{\dot x}^2 +
2\Gamma^2_{12}{\dot x}{\dot y} + \Gamma^2_{22}{\dot y}^2\bigr), \\
\end{aligned}
\right.
\label{4}
\end{equation}
where $\Del(x,y) = ac-b^2$ is the discriminant of the metric and $\Gamma^k_{ij}$ are the Christoffel symbols.

Clearly, the definition of geodesics as auto-parallel curves in the Levi-Civita connection
generated by the metric \eqref{1} leads to the same Equation~\eqref{4}, see e.g.~\cite{DNF}.

\medskip

We shall distinguish three possible types of geodesics and, more generally, three types of curves in the given metric:

\begin{defin}
A curve $\gam$ is called timelike (spacelike) if the inequality $ds^2 > 0$
($ds^2 <0$, respectively) holds true at all points of $\gam$.
A curve $\gam$ is called isotropic if $ds^2 \equiv 0$ along $\gam$.
\end{defin}


Equation \eqref{3}, or equivalently, \eqref{4} defines geodesics as parametrized
curves (the corresponding parametrization is called natural), which are images of integral curves of
the field
\begin{equation}
2\Del \biggl( \dot x \frac{\pa}{\pa x} +  \dot y \frac{\pa}{\pa y}
\biggr) + (cP-bR) \frac{\pa}{\pa \dot x} + (aR-bP) \frac{\pa}{\pa \dot y}
\label{5}
\end{equation}
under the projection $\pi_1 \colon T\M \to \M$.

\medskip

In a domain where the metric \eqref{1} is smooth and non-degenerate
one can exclude from consideration the zero section of the tangent bundle $T\M$,
since it consists of equilibrium points of \eqref{4}.
Hence at every non-parabolic point all geodesics have definite tangential directions
$p = \frac{\dot y}{\dot x} = \frac{dy}{dx}$.
However, in general it is impossible to assert the same for parabolic points.
Further we restrict ourselves to geodesics with definite tangent directions at parabolic points.

In view of what we have said above, unparametrized geodesics can be (locally) defined as extremals of the length functional $\JL(\gamma) = \int_{\gamma} \sqrt{a + 2bp + cp^2}\,dx$.
The corresponding Euler-Lagrange equation reads
\begin{equation}
\frac{d}{dx} \LL_{p} - \LL_{y} = 0,
\quad \textrm{where} \quad
\LL(x,y,p) = \sqrt{F}, \ \, F = a + 2bp + cp^2.
\label{6}
\end{equation}
The Lagrangian $\LL(x,y,p)$ is a function of the projective tangent bundle $PT\M$,
and Equation \eqref{6} is correct at all points of $PT\M$ except the isotropic surface $\FF\colon F(x,y,p)=0$.

Straightforward transformations of Equation \eqref{6} show that unparametrized geodesics are images of integral curves of the field $\W$ given by the formula
\begin{equation}
\W = \frac{1}{2 F^{\frac{3}{2}}} \, \V,
\quad \textrm{where} \quad
\V = 2\Del \biggl(\frac{\pa}{\pa x} + p \frac{\pa}{\pa y} \biggr) + M \frac{\pa}{\pa p},
\ \,
M = \sum\limits_{i=0}^{3} \mu_i(x,y)p^i,
\label{7}
\end{equation}
\begin{equation}
\begin{aligned}
\mu_3 &= c(2b_y-c_x) - bc_y, \ \
\mu_2 = b(2b_y-3c_x) + 2a_yc - ac_y, \\
\mu_1 &= b(3a_y-2b_x) + a_xc - 2ac_x, \ \ \mu_0 = a(a_y-2b_x) +a_xb,
\end{aligned}
\label{8}
\end{equation}
under the projection $\pi_2 \colon PT\M \to \M$.

The spaces $T\M$ and $PT\M$ are connected with the mapping
$\Pi \colon T\M \to PT\M$, the projectivization of tangent planes to the
manifold $\M$. The relationship between the projections $\pi_1$ and
$\pi_2$ is presented on the scheme:
\begin{equation*}
\begin{matrix}
\xymatrix{ {T\M}\ar[dr]^-{\displaystyle \pi_1} \ar[dd]_{
{{\displaystyle \Pi}} \phantom{.}
}& \\
& {\M}\\
{PT\M}\ar[ur]_{\displaystyle \pi_2} &
}\\
\end{matrix}
\end{equation*}


\begin{theorem}
\label{T1}
$\Pi \colon T\M \to PT\M$ sends the field \eqref{5} to the field $\V$ given by the formula \eqref{7}.
\end{theorem}

{\bf Proof}. The mapping $\Pi$ assigns to a tangent vector $(\dot x,
\dot y)$ the value $p={\dot y}/{\dot x}$, therefore
\begin{equation}
\frac{dp}{dx} = \frac{dp}{dt} \biggl( \frac{dx}{dt}\biggr)^{-1} =
\frac{1}{\dot x}\, \frac{d}{dt} \biggl( \frac{\dot y}{\dot x}
\biggr) = \frac{1}{{\dot x}^3}\, \biggl( \dot x \frac{d\dot y}{dt} -
\dot y \frac{d\dot x}{dt}   \biggr) = \frac{1}{{\dot x}^2}\, \biggl(
\frac{d\dot y}{dt} - p \frac{d\dot x}{dt} \biggr).
\label{9}
\end{equation}
From \eqref{5} we obtain the relations
\begin{equation*}
\frac{d\dot x}{dt} = \frac{cP-bR}{2\Del}, \quad  \frac{d\dot y}{dt}
= \frac{aR-bP}{2\Del},
\end{equation*}
and from \eqref{3} we have
\begin{equation*}
\frac{P}{{\dot x}^2} = (c_x-2b_y) p^2 - 2a_y p - a_x, \ \
\frac{R}{{\dot x}^2} = -c_y p^2 - 2c_x p + a_y - 2b_x.
\end{equation*}
Substituting the above expressions in \eqref{9}, after straightforward
transformations we obtain
$$
\frac{dp}{dx} = \frac{1}{2\Del} \biggl( \frac{aR-bP}{{\dot x}^2} -
p\,\frac{cP-bR}{{\dot x}^2} \biggr) = \frac{1}{2\Del}
\sum\limits_{i=0}^{3} \mu_i p^i,
$$
where the coefficients $\mu_i$ are defined in \eqref{8}.
\hfill\ebox

\medskip

For a Riemannian metric, Equation \eqref{6} is well defined on the whole $PT\M$, and naturally parametrized geodesics (i.e., extremals of the action functional) can be obtained from extremals of the length
functional by means of the arc-length parametrization.
For a Lorenzian metric, Equation \eqref{6} is not defined on the isotropic surface $\FF$ and
the arc-length parametrization for isotropic curves is impossible.
However, as we shall see in the sequel, the notion of natural parametrization
in Lorenzian metrics does make sense for all geodesics, including isotropic ones.

\begin{theorem}
\label{T2}
The isotropic surface $\FF$ is an invariant surface of the field $\V$.
Isotropic curves with an appropriate choice of parametrization (called natural)
are extremals of the action functional, that is, they are geodesics.
\end{theorem}

{\bf Proof}. Since the vector field $\W$ from the formula \eqref{7}
is obtained directly from the Euler-Lagrange equation \eqref{6},
the divergence of $\W$ is identically equal to zero everywhere except for the isotropic
surface $\FF$. By~\cite{GR} (Theorem~1), this implies that $\FF$ is an invariant
surface of the field~$\V$.

On the isotropic surface $\FF$, contact planes $dy=pdx$ cut a
direction field $\mathcal X$, whose integral curves are 1-graphs of
isotropic curves in the given metric. On the other hand, $\mathcal X$ coincides with the
restriction of the field $\V$ to $\FF$, hence the projection $\pi_2
\colon PT\M \to \M$ sends integral curves of the field $\mathcal X$ to (unparametrized) geodesics.
\hfill\ebox

\begin{remark}
\label{R1}
The first assertion of Theorem~\ref{T2} is valid for any $\dim \M > 2$,
while the second assertion is valid only for $\dim \M = 2$,
see the following example.
\end{remark}


{\bf Example 1.1.}
Consider the metric $ds^2 = dx^2 + dy^2 - dz^2$. Geodesics are various straight lines, and
isotropic geodesics are generatrices of the cones $(x-x_0)^2 + (y-y_0)^2 - (z-z_0)^2 = 0$.
However, there exist isotropic curves that are not geodesics, for instance, the
elliptic helix $x=\sin t$, $y=\cos t$, $z=t$ or the hyperbolic helix
$x=\sh t$, $y=t$, $z=\ch t$.
\hfill\ebox

\begin{remark}
\label{R2}
From Theorem~\ref{T2} it is seen that geodesics keep their types
(remain timelike, spacelike or isotropic)
on any interval that contains no parabolic points.
\end{remark}



\section{Local properties of geodesics}\label{Local}

\subsection{Unparametrized geodesics}\label{Unparametrized}

In this section, we briefly describe the local properties of unparametrized geodesics
in pseudo-Riemannian metrics in a neighborhood of generic parabolic points
(for proofs, see \cite{GR, Rem-Pseudo}).
More precisely, here we use an auxiliary parametrization of geodesics,
which is not the natural parametrization defined by Equation~\eqref{2}.

Consider the metric \eqref{1} with smooth ($C^{\infty}$) coefficients $a,b,c$ that is Riemannian
on an open domain, Lorenzian on some other open domain, and degenerate on the curve
$\AA \colon \Del(x,y)=0$, which separates these domains.
In following the paper \cite{Mier}, we give the following definition:

\begin{defin}
A parabolic point $q_0 \in \AA$ is called transverse if at this point
the coefficients $a,b,c$ do not vanish simultaneously,
the discriminant $\Del$ has nonzero gradient, and
the unique isotropic direction $p_0 = -\frac{a}{b}$ is transversal to the degeneracy curve $\AA$,
that is, $b\Del_x - a\Del_y \neq 0$.
\end{defin}

For any transverse parabolic point $q_0 = (x_0,y_0)$, the unique isotropic direction $p_0$ is a prime root
of the cubic polynomial $M(q_0,p)$ with respect to the variable $p \in \PR$.
Generically, the cubic polynomial $M(q_0,p)$ has one or three prime roots: $p_i \in \PR$, $i=0$ or $i=0,1,2$.

\begin{theorem}
\label{T3}
For a generic transverse parabolic point $q_0$, there exist one or three admissible directions,
which correspond to roots of the cubic polynomial $M(q_0,p)$, $p \in \PR$.
To the isotropic direction $p_0$ would correspond an infinite number of geodesics indexed by a real parameter.
To each non-isotropic admissible direction $p_i$, $i=1,2$, would correspond a unique geodesic.
Almost all geodesics with the isotropic tangential direction have a cusp at $q_0$, while geodesics with non-isotropic tangential directions are regular.
\end{theorem}

\begin{remark}
\label{R3}
Let \eqref{1} be a metric induced on a smooth surface $\M$ embedded in the three-dimensional Minkowski space.
Then there exists a very graphic interpretation of the cases when the cubic polynomial $M(q_0,p)$ has one or three real roots: Gaussian curvature of $\M$ at the point $q_0$ is positive or negative, respectively.
\end{remark}

\medskip

\begin{figure}[h!]
\begin{center}
\includegraphics{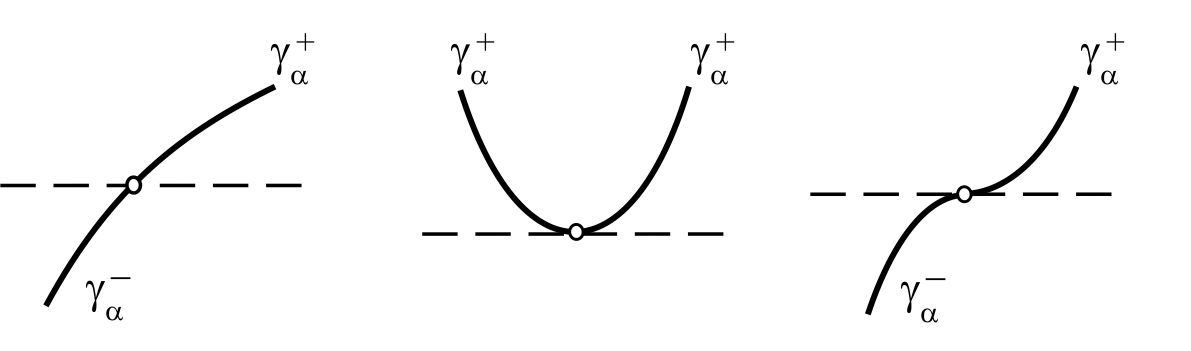}
\caption{
A geodesic $\gama$ passing through a given point $q_0$ with a non-horizontal tangential direction
corresponds to the couple $(\gamap, \gamam)$ (left). A~geodesic $\gama$ with the horizontal direction at $q_0$
corresponds to the couple $(\gamap, \gamap)$ (center) or $(\gamam, \gamam)$ or $(\gamap, \gamam)$ (right).}
\label{Fig1}
\end{center}
\end{figure}

Further we distinguish geodesics outgoing from a point (not necessary parabolic)
in the semiplanes $y>y_0$ and $y<y_0$ using the superscripts $\pm$ (see Fig.~\ref{Fig1}).

Let $q_0 \in \M$ be a transverse parabolic point.
Without loss of generality assume that $q_0=0$ (the origin) and the conditions $a(0) > 0$, $b(0)=c(0)=0$ hold.
Then the isotropic direction $p_0=\infty$ and the condition $b\Del_x - a\Del_y \neq 0$ reads $c_y(0) \neq 0$.
For definiteness, suppose that $c_y(0) < 0$.

By $\Gamo$ denote the family of geodesics outgoing from the point $q_0$ with the isotropic tangential direction
indexed by a certain parameter $\alf$, that is, $\Gamo = \{\gamapm\}$.


\begin{theorem}
\label{T4}
There exist smooth local coordinates centered at $0$ such that:

1. \
The degeneracy curve $\AA$ coincides with the $x$-axis,
the semiplane $y>0$ is Lorenzian and the semiplane $y<0$ is Riemannian.

2. \
Geodesics $\gamma^{+}_{\pm \alf}$ $(\gamma^{-}_{\pm \alf})$
are the branches of the semicubic parabolas
\begin{equation}
x= \alf \tau^3 X_{\alf}(\tau), \quad y= \tau^2 Y_{\alf}(\tau), \quad
\tau \in (\bR, 0)_{\alf}, \quad \alf \ge 0,
\label{10}
\end{equation}
where $X_{\alf}, Y_{\alf}$ are smooth functions, $X_{\alf}(0)=1$, $Y_{\alf}(0)=+1$ $(-1)$,
and $(\bR, 0)_{\alf}$ are neighborhoods of zero depending on $\alf$.
Here the superscript $\pm$ of a geodesic coincides with the sign of $Y_{\alf}(0)$,
while $\pm$ in the subscript corresponds to the left $(\tau < 0)$ and right $(\tau > 0)$ branches of the semicubic parabola.\footnote{
Here the case $\alf = 0$ is interpreted as the limiting case of semicubic parabolas
(the branches are glued together), and geodesics $\gam^{\pm}_0$ are the halves of the $y$-axis.
}
\end{theorem}

\begin{figure}[h!]
\begin{center}
\includegraphics{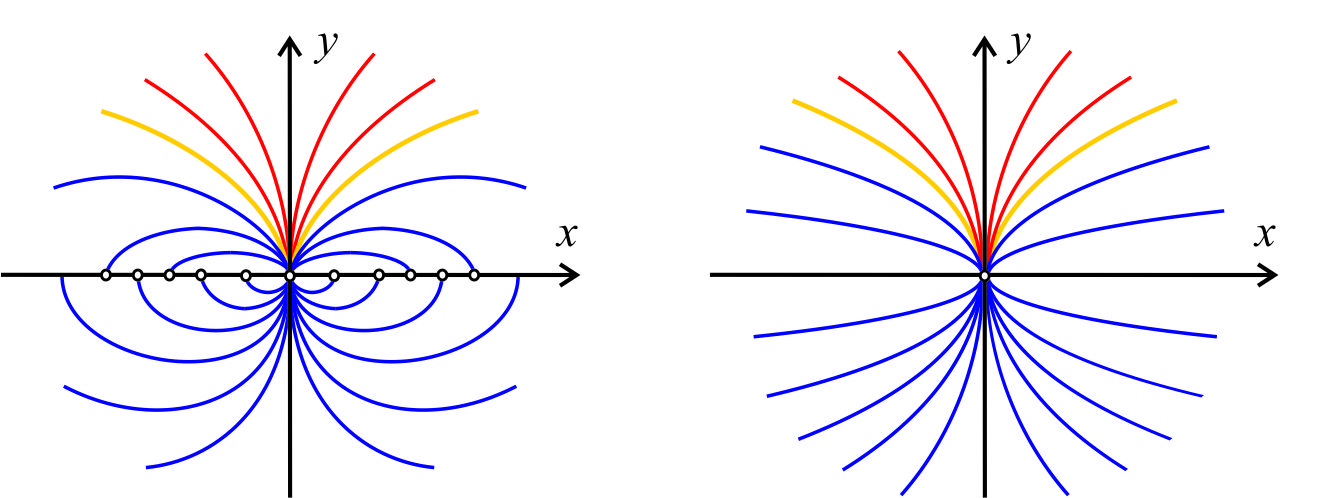}
\caption{
Two examples of the phase portraits of geodesics outgoing from the origin: timelike (blue lines), spacelike (red lines),
and isotropic (yellow lines).
}
\label{Fig5}
\end{center}
\end{figure}

Formula \eqref{10} shows that the family $\Gamo$ contains geodesics of all three possible types:
timelike ($\gamap$ with $|\alf| > \alfis$ and all $\gamam$),
spacelike ($\gamap$ with $|\alf| < \alfis$),
and isotropic ($\gamap$ with $|\alf| = \alfis$),
where $\alfis = \tfrac{2}{3} \sqrt{-c_y(0)/a(0)}$.
The possible mutual relationship between them is presented on Fig.~\ref{Fig5}.

{\bf Example 2.1.}
An example of the phase portrait presented on Fig.~\ref{Fig5} (right) is given by the metric
$ds^2 = dx^2 + y dy^2$. Geodesics of the family $\Gamo$ are defined by the formula $x = \alf \tau^3$, $y = \pm \tau^2$ for all $\tau \in \bR$. Not a single geodesic will return to $\AA$ (the $x$-axis).
An example of the phase portrait presented on Fig.~\ref{Fig5} (left) is given by the metric
$$
ds^2 = \frac{1}{1+y^2}\,dx^2 + y dy^2.
$$
Geodesics in this metric can be defined from formula~\eqref{6}.
One can see that all geodesics of the family $\Gamo$ (except only vertical straight lines) return to $\AA$ (the $x$-axis). We omit the proof of this statement, because it is a trivial corollary of a more general fact (Theorem~\ref{T6}).
\hfill\ebox

This leads to the following crucial question:
\textit{
Does a given geodesic $\gamapm$ outgoing from a transverse parabolic point $q_0 \in \AA$
return to the degeneracy curve $\AA$, or not?
}
It is not easy (if even possible) to give the complete answer for arbitrary pseudo-Riemannian metrics.
In Section~\ref{Pseudo-Riemannian}, we shall give the answer for metrics possessing differentiable groups of symmetries, for instance, if the coefficients do not depend on one of the coordinates.

\medskip





\subsection{Naturally parametrized geodesics}\label{Parametrized}

In this section, we consider naturally parametrized geodesics, that is, consider them
as extremals of the action functional (Equation~\eqref{3}).
In particular, we prove that all geodesics $\gamapm \in \Gamo$ reach the parabolic point $q_0$ in finite time, assuming that the natural parametrization of $\gamapm$ is chosen so that the motion proceeds toward~$q_0$.

Note that the parametrization \eqref{10} from Theorem~\ref{T4} is not unique:
there is an infinite group of reparametrizations $\tau \to \tau \phi(\tau)$ with arbitrary smooth functions $\phi$ that $\phi(0)=1$. Theorem~\ref{T5} below asserts that among all possible parametrizations \eqref{10}
there exists one (and only one) that is connected with the natural parametrization
through the simple relation $\tau = t^{\frac{1}{3}}$.

\begin{theorem}
\label{T5}
All geodesics $\gamapm \in \Gamo$ reach a transverse parabolic point $q_0$ in finite time.
In the local coordinates from Theorem~\ref{T4}, the natural parametrization of $\gamma^{+}_{\pm \alf}$ $(\gamma^{-}_{\pm \alf})$ has the form
\begin{equation}
x= \alf t X_{\alf}(t^{\frac{1}{3}}),
\quad
y= t^{\frac{2}{3}} Y_{\alf}(t^{\frac{1}{3}}),
\quad t \in (\bR, 0), \quad \alf \ge 0,
\label{13}
\end{equation}
where $X, Y$ are smooth functions, $X(0)=1$, $Y(0)=+1$ $(-1)$, up to scalings of $t$.
\end{theorem}

{\bf Proof}.
Assume $q_0=0$ and consider a couple $\gamma^{+}_{\pm \alf}$ with fixed $\alf > 0$, consisting a sole semicubic parabola from \eqref{10}. (For $\gamma^{-}_{\pm \alf}$ the proof is similar.)
Then $Y_{\alf}(0)= 1$, and after an appropriate reparametrization,
the semicubic parabola has the form
$x= \tau^3 (\alf + f_{\alf}(\tau))$, $y= \tau^2$,
with a smooth function $f_{\alf}$ vanishing at zero.
Represent the germ of $f_{\alf}$ in the form
$f_{\alf}(\tau) = \phi_{\alf}(\tau^2) + \tau \psi_{\alf}(\tau^2)$
with smooth functions $\phi_{\alf}, \psi_{\alf}$.
In a neighborhood of the origin, the change of variables
\begin{equation}
x \to \frac{x - y^2 \psi_{\alf}(y)}{\alf + \phi_{\alf}(y)},  \quad \phi_{\alf}(0)=0,
\label{11}
\end{equation}
brings the semicubic parabola to the form $x= \tau^3$, $y= \tau^2$, or equivalently, $y= x^{\frac{2}{3}}$.

Without loss of generality, assume $a(0)=1$. Then, substituting
$$
y= x^{\frac{2}{3}}, \quad
\dot y= \frac{2}{3} x^{-\frac{1}{3}} \dot x, \quad
\ddot y= \frac{2}{3} x^{-\frac{1}{3}} \ddot x - \frac{2}{9} x^{-\frac{4}{3}} {\dot x}^2
$$
in the first equation of \eqref{3}, after transformations we get the equation
\begin{equation}
\frac{\ddot x}{\dot x} =
\frac{a_1 x^{-\frac{1}{3}} + a_2 x^{-\frac{2}{3}} + \psi(x^{\frac{1}{3}})}{1 + x^{\frac{1}{3}} \phi(x^{\frac{1}{3}})}
\, \dot x,
\label{14}
\end{equation}
where $a_1, a_2$ are certain constants and $\phi, \psi$ are smooth functions.
Integrating \eqref{14}, we get
\begin{equation}
\ln{|\dot x|} = F(x) + \const, \quad
F(x) = \int\limits_{0}^{x} \frac{a_1 \xi^{-\frac{1}{3}} + a_2 \xi^{-\frac{2}{3}} + \psi(\xi^{\frac{1}{3}})}{1 + \xi^{\frac{1}{3}} \phi(\xi^{\frac{1}{3}})} \, d\xi = f(x^{\frac{1}{3}}),
\label{145}
\end{equation}
where $f$ is a smooth function, $f(0)=0$.

Equation \eqref{145} yields $\dot x = C e^{F(x)}$.
Without loss of generality we put $C = 1$, this corresponds
to a choice of the velocity of the motion along the geodesic.
Then we arrive at the differential equation
$$
\frac{dt}{dx} = e^{-F(x)} = 1 - F(x) + \frac{1}{2}F^2(x) + \cdots =
1 + g(x^{\frac{1}{3}}),
$$
which has the solution
\begin{equation}
t = x (1 + h(x^{\frac{1}{3}})) + \const,
\label{15}
\end{equation}
where $g, h$ are smooth functions, $g(0)=h(0)=0$.
The equality \eqref{15} shows that both geodesics $\gamma^{+}_{\pm \alf}$
reach the parabolic point $0$ in finite time (this does not depend on the choice of local coordinates).

To prove the second statement of the theorem, we put $\const=0$ in \eqref{15} and express
the variable $x$ through $t$. Then we get
\begin{equation}
x = t \wti X(t^{\frac{1}{3}}), \quad y = t^{\frac{2}{3}} \wti Y(t^{\frac{1}{3}})
\label{155}
\end{equation}
with certain smooth functions $\wti X, \wti Y$, $\wti X(0) = \wti Y(0) =1$.
The expression \eqref{155} is established in the local coordinates given by~\eqref{11}.
Apply the inverse transformation $x \to (\alf + \phi_{\alf}(y))x + y^2 \psi_{\alf}(y)$,
we obtain~\eqref{13} in the initial coordinates.

There remains one case to consider, that in which $\alf = 0$, and both geodesics $\gamma^{+}_{\pm \alf}$
coincide with the positive half of the $y$-axis.
Substituting $x=0$ in the second equation of \eqref{3}, we get
\begin{equation}
2 \, \frac{\ddot y}{\dot y} = - \frac{\phi'(y)}{\phi(y)} \, \dot y,
\label{16}
\end{equation}
where $\phi$ is a smooth function, $\phi(0)=0$, $\phi'(0)\neq 0$.
Using analogous reasonings for \eqref{16} as were used for \eqref{14},
one can prove the theorem in the case $\alf = 0$.
\hfill\ebox

{\bf Example 2.2.}
Consider geodesics in the metric $ds^2 = dx^2 - y dy^2$. Then the system \eqref{3} consists of two equations
$\ddot x = 0$, $2y \ddot y + {\dot y}^2 = 0$,
and the integration yields $x = c_1 t$, $y = c_2 t^{\frac{2}{3}}$.
Hence naturally parametrized geodesics $\gamma^{+}_{\pm \alf}$ $(\gamma^{-}_{\pm \alf})$
have the form $x = \alf t$, $y = t^{\frac{2}{3}}$ $(-t^{\frac{2}{3}})$, $\alf \ge 0$,
up to scalings of $t$.
We thus see again that the natural parametrization (due to the action functional) is defined for all geodesics including isotropic, while the arc-length parametrization for them does not exist.
\hfill\ebox


\section{Metrics with differentiable groups of symmetries}\label{Symmetries}

In this section, we study global properties of geodesics in metrics
\begin{equation}
ds^2 = a(y) \, dx^2 + 2b(y) \, dx dy + c(y) \, dy^2,
\label{18}
\end{equation}
conditioned by the existence of the group of symmetries consisting of parallel translations along the $x$-axis.

Let $q_0 = (x_0,y_0)$ be an arbitrary point and $\Gamo = \{\gamapm\}$ be the family of geodesics
outgoing from $q_0$ indexed by a parameter $\alf \in \bR$.
(As before, we distinguish geodesics outgoing in the semiplanes $y>y_0$ and $y<y_0$ using the superscripts $\pm$.)
We shall examine geodesics of the family $\Gamo$ in an open domain $\ov \Om$
defined by the inequalities $\om^- <y< \om^+$, where $\om^{\pm}$ can be finite or infinite.\footnote{
In particular, $\ov \Om$ can coincide with the whole $(x,y)$-plane: $\om^{\pm} = \pm \infty$.
}
By $\Om$ denote the domain obtained from $\ov \Om$ by eliminating the line $y = y_0$, where $\om^- <y_0< \om^+$.

\begin{assum}
\label{AS1}
In the domain $\Om$ the coefficient $a(y)$ does not vanish and the metric \eqref{18} is smooth and nondegenerate
(not necessarily pseudo-Riemannian). However, these conditions can be not satisfied on the boundary of $\Om$, for instance, on the line $y=y_0$.
\end{assum}


\begin{lemma}
\label{L2}
Every unparametrized geodesic $\gam \subset \Om$ is a solution of the equation
\begin{equation}
(b^2-h^2c) p^2 + 2b(a-h^2) p + a(a-h^2) = 0,  \quad p = \frac{dy}{dx},
\label{19}
\end{equation}
with appropriate constant $h^2$, where $h^2 \ge 0$ ($h^2 \le 0$) if $\gam$ is timelike (respectively, spacelike) and
$h^2=\infty$ if $\gam$ is isotropic. In the case $h^2=\infty$ the equation \eqref{19} reads $cp^2+ 2bp +a = 0$.
\end{lemma}

{\bf Proof}.
By Theorem~\ref{T2} (see also Remark~\ref{R2}), the function $F = a + 2bp + cp^2$ keeps its sign along every geodesic $\gam \subset \Om$. Then the Lagrangian $\LL = \sqrt{F}$ does not vanish or is identically equal to zero along $\gam$. In the last case we have the equation of isotropic geodesics $cp^2+ 2bp +a = 0$.

Now assume that $\LL$ does not vanish along $\gam \subset \Om$.
The equation of unparametrized geodesics is the Euler--Lagrange equation \eqref{6}.
Since the Lagrangian $\LL$ does not depend on $x$, the equation \eqref{6} has the energy integral
\begin{equation}
H(y,p) = \LL - p \LL_p = \frac{a+bp}{\sqrt{a+2bp+cp^2}},
\label{20}
\end{equation}
and $\gam$ satisfies the equation $H(y,p) = h$, where $h$ is an appropriate constant.
After straightforward transformations, this yields \eqref{19}.
\hfill\ebox

However, not every solution of the equation \eqref{19} is a solution~\eqref{6}.

{\bf Example 3.1.}
The Klein metric
\begin{equation}
ds^2 = \frac{dx^2 + dy^2}{y^2}
\label{21}
\end{equation}
satisfies all above conditions in $\Om$, which is the $(x,y)$-plane without the line $y=0$.
Geodesics in this metric are the circles $(x-x_*)^2 + y^2 = r^2$, $r \neq 0$, and vertical straight lines $x=x_*$.
For the metric \eqref{21}, the equation \eqref{19} reads
\begin{equation}
h^2 p^2 + h^2 - y^{-2} = 0,  \quad p = \frac{dy}{dx}.
\label{2100}
\end{equation}
The value $h^2=0$ gives the lines $x=x_*$, while $h^2>0$ gives the circles $(x-x_*)^2 + y^2 = r^2$ with $r=h^{-1}$.
It is not hard to see that for any $h^2>0$ the equation \ref{2100} has also the constant solution $y \equiv h^{-1}$,
which is an envelop of the circles $(x-x_*)^2 + y^2 = h^{-2}$, but not a solution of \eqref{6}, and consequently, not a geodesic in the metric \eqref{21}.
\hfill\ebox

\subsection{Implicit differential equations}\label{Symmetries-IDE}

Equation \eqref{19} with a fixed constant $h$ is a partial case of so-called {\it implicit} differential equation, that  is, a differential equation not solved with respect to the derivative:
\begin{equation}
\F(x,y,p) = 0, \quad p = \frac{dy}{dx}.
\label{22}
\end{equation}
The equation \eqref{22} defines a multi-valued direction field on the $(x,y)$-plane, which becomes
single-valued on the surface $\F(x,y,p)=0$ in the $(x,y,p)$-space.
The passage from the multi-valued field on the plane to the single-valued {\it lifted field} on the surface is similar to the Riemann surface for a multi-valued function of a complex variable on which this function becomes single-valued.

The lifted field is an intersection of the contact planes $dy = pdx$ with the tangent planes to the surface $\F(x,y,p)=0$, that is, it is defined by the vector field
\begin{equation}
\F_p\frac{\pa}{\pa x} + p\F_p\frac{\pa}{\pa y} - (\F_x + p\F_y)\frac{\pa}{\pa p}.
\label{220}
\end{equation}
The lifted field is defined at all points where the surface $\F(x,y,p)=0$ is regular and the contact plane is not tangent to this surface. Solutions of the equation \eqref{22} are curves whose 1-graphs are integral curves of the field \eqref{220}.

Points of the surface $\F(x,y,p)=0$ where $\F_p = 0$ are called {\it singular},
the set of singular points is called the {\it criminant}, and the projection of the criminant on the $(x,y)$-plane is called the {\it discriminant curve}. A singular point is called {\it proper} if $\F_x + p \F_y \neq 0$, otherwise, it is called {\it improper}.
Two implicit differential equations are called smoothly (topologically) {\it equivalent} if there
exists a diffeomorphism (homeomorphism, respectively) of the $(x,y)$-plane that sends integral
curves of the first equation to integral curves of the second one.
The smooth (topological) local classification of implicit differential equations in a neighborhood of their singular points is an important problem of the qualitative theory of differential equations; for more retails see  \cite{Dav-Japan, Rem2006}.

For a generic implicit differential equation \eqref{22} almost all singular points are proper, and
improper points (being saddles, nodes or foci of the lifted field) are isolated on the criminant.
However, the equation \eqref{19}, which is of interest to us, is not generic, since its left-hand side
does not depend on $x$. A (non-generic) implicit differential equation may have so-called {\it singular solutions} being envelops of the family of all solutions. It is not hard to see that the 1-graph of a singular solution consist of improper singular points. Consequently, any singular solution is contained in the discriminant curve of the equation. For example, the implicit differential equation \eqref{2100} with $h \neq 0$ has the singular solution $y \equiv h^{-1}$, which is not a geodesic in the metric \eqref{21}. Here we become acquainted with the general fact: as we shall soon see, solutions of the equation \eqref{19} not being geodesics in the metric \eqref{18} are exactly singular solutions of \eqref{19}.

\subsection{Basic results}\label{Symmetries-Basic}

In the following two lemmas, we consider the equation \eqref{19} with arbitrary fixed $h \neq 0, \infty$, in the domain $\Om$. By $\F$ denote the left-hand side of this equation.

\begin{lemma}
\label{L3}
The discriminant curve is given by the equation $a(y)=h^2$, the criminant is the intersection of the surface  $\F(x,y,p)=0$ and the plane $p=0$. At any singular point, $\F_{pp} \neq 0$ holds true and $\F_{y} \neq 0$ holds if and only if $a'(y) \neq 0$.
\end{lemma}

{\bf Proof}.
Consider the equation \eqref{19} as a quadratic equation in $p$ with the discriminant
\begin{equation}
D = (a-h^2) (ac-b^2) h^2
\label{23}
\end{equation}
depending on $y$.
Since $h \neq 0$ and $\Del(y) \neq 0$ in $\Om$, the discriminant curve is given by the equation $a(y)=h^2$.
Hence at each point of the discriminant curve $b^2-h^2c = - \Del \neq 0$, $\F_{pp} \neq 0$,
and the equation \eqref{19} has the double root $p=0$.
Finally, on the the discriminant curve we have the equalities $\F_{y}(y,0) = a(y)a'(y)$ and $a(y)=h^2 \neq 0$.
\hfill\ebox

Lemma \ref{L3} asserts that the equation \eqref{19} may have singular solutions only in the form $y = \const$.
On the other hand, any $y = \const$ is a solution of the equation \eqref{19} with an appropriate $h^2 \neq \infty$
(substituting $y \equiv y_*$ in \eqref{19}, it is easy to see that the equation is satisfied with $h^2 = a(y_*)$).
Example~3.1 shows that a singular solution of \eqref{19} may be not a geodesic.
The following lemma gives a criterion whether $y = \const$ is a geodesic in the metric \eqref{18}.

\begin{lemma}
\label{L4}
Let $a(y_*)=h^2$, then there is the following alternative:
Either $a'(y_*) = 0$, the constant function $y \equiv y_*$ is a non-singular solution of \eqref{19} and it is a non-isotropic geodesic, or else $a'(y_*) \neq 0$, the function $y \equiv y_*$ is a singular solution of \eqref{19} and it is not a geodesic. In the second case, the equation \eqref{19} is smoothly equivalent to ${\ti p}^2 = {\ti y}$ in a neighborhood of any point on $y=y_*$.
\end{lemma}

{\bf Proof}.
Substituting $y \equiv y_*$ in the equations \eqref{19} and \eqref{6}, one can see that $y \equiv y_*$ is a solution of \eqref{19} and it is a geodesic if and only if $a'(y_*) = 0$ (since the constant $h^2 = a(y_*) \neq \infty$, it is non-isotropic). Thus it remains to prove that $y \equiv y_*$ is a singular solution of \eqref{19} if and only if $a'(y_*) \neq 0$.

Let $a'(y_*) \neq 0$. By Lemma~\ref{L3}, then $F(0,y_*) = F_{p}(0,y_*) = 0$ and $F_{pp}(0,y_*) \neq 0$, $F_{y}(0,y_*) \neq 0$. A germ of implicit differential equations satisfying these conditions is called {\it Clairaut fold}, it is smoothly equivalent to ${\ti p}^2 = {\ti y}$; see~\cite{Dav-Japan}.
Here the discriminant curve $y=y_*$ of the equation \eqref{19} corresponds to the discriminant curve ${\ti y}=0$ of the normal form ${\ti p}^2 = {\ti y}$. The last one is a singular solution: the line ${\ti y}=0$ is an envelop of  parabolas ${\ti y}= \frac{1}{4} ({\ti x}-c)^2$, $c \in \bR$, being non-singular solutions of ${\ti p}^2 = {\ti y}$.

Now let $a'(y_*) = 0$. As we have shown, $y \equiv y_*$ is a geodesic and it is a solution of the equation \eqref{19} with $h^2 = a(y_*)$. Assume that $y \equiv y_*$ is a singular solution, i.e., it is an envelop of non-singular solutions of \eqref{19}, which therefore are not geodesics. Indeed, otherwise there exist two different geodesics passing through each point of the line $y \equiv y_*$ with the same tangential direction. This is impossible, since he line $y \equiv y_*$ belongs to the domain $\Om$, where the metric is smooth and non-degenerate.
On the other hand, if non-singular solutions of \eqref{19} are not geodesics, there exist points $(x,y) \in \Om$ and tangential directions $p \in \PR$ such that there is no geodesics passing through $(x,y)$ with given tangential direction $p$. The contradiction proves that $y \equiv y_*$ is a non-singular solution of the equation~\eqref{19}.
\hfill\ebox

\medskip

Without loss of generality, assume that the family $\Gamo$ is indexed by $\alf$ so that both geodesics $\gamapm \in \Gamo$ correspond to the same value $\ha^2$, i.e., the equality \eqref{19} with $h^2=\ha^2$ holds identically on both $\gamapm$. It is always possible in all cases of interest. For instance, if $q_0$ is a transverse parabolic point of a pseudo-Riemannian metric, the parameter $\alf$ is chosen as in Theorems~\ref{T4},~\ref{T5}.
If the metric is smooth and non-degenerate at $q_0$, one can put $\alf$ equal to the value $p$ or $p^{-1}$ at $q_0$.

For $p=0$, the equality \eqref{19} gives $a(y)=\ha^2$.
Consider the least equation with respect to the unknown $y$. Put
\begin{equation*}
\yap = \inf \{ y_0<y<\om^+ \colon a(y)=\ha^2 \},
\ \ \,
\yam = \sup \{ \om^-<y<y_0 \colon a(y)=\ha^2 \},
\end{equation*}
and $\yapm = \om^{\pm}$ if the corresponding set is empty.

\begin{assum}
\label{AS2}
We shall assume further that $\yam < y_0 < \yap$ and there are no geodesics asymptotically tending to the horizontal line $y=y_0$. The last condition holds true in all cases of interest, including Riemannian and Lorentzian metrics
as well as pseudo-Riemannian metrics with transverse parabolic points.
\end{assum}

\begin{lemma}
\label{L5}
Let $y^{\pm}$ be such that $y_0 < y^+ < \yap$ and $\yam < y^- < y_0$.
Then any geodesic $\gamapm \in \Gamo$ with $\ha^2 \neq 0, \infty$ reaches the horizontal line $y=y^{\pm}$ (with the same superscript) at a finite point.
\end{lemma}

{\bf Proof}.
For definiteness, consider $\gamap$ (for $\gamam$ the proof is similar).
Since $\gamap$ goes out from the point $q_0$ to the semiplane $y>y_0$,
there exists $\ye > y_0$ such that $\gamap$ intersects the line $y=\ye$.
If $\ye \ge y^+$, then the required statement is proved.

Suppose that $\ye < y^+$.
Then $a(y) \neq 0$, $a(y)-\ha^2 \neq 0$ and $D(y)>0$ for all $y \in Y = [\ye, y^+]$.
Hence for any $y \in Y$, $p=0$ is not a root of the quadratic equation \eqref{19},
and the function $P(y) = \inf \, \{ |p| \colon H^2(y,p) = \ha^2 \} > 0$ for all $y \in Y$.
$P(y)$ reaches a minimum value $p_* > 0$ on $Y$, and the geodesic $\gamap$ reaches the horizontal line $y=y^+$ not later that the length of its arc (beginning at the intersection with the line $y=\ye$ and considered in the standard Euclidean metric $dx^2 + dy^2$) reaches the value $(y^+ - \ye)/p_*$.
\hfill\ebox

Before formulating the following theorem, we make an important convention.
Assume that a geodesic $\gam$ leaves the domain $\Om$ through the boundary $y = \om^{\pm}$ and then
it returns again in $\Om$. By $\gam'$ and $\gam''$ denote the corresponding parts of $\gam$ lying in $\Om$.
Then we shall consider $\gam'$ and $\gam''$ to be different geodesics.
So, if a geodesic intersects the line $y = \om^{\pm}$, it never returns back in $\Om$.

\begin{theorem}
\label{T6}
A geodesic $\gamap \in \Gamo$ with any $\ha^2$ returns to the initial horizontal line $y=y_0$
if and only if $\yap < \om^+$ and $a'(\yap) \neq 0$.
Moreover, the following triple choice holds:

1. If  $\yap < \om^+$ and $a'(\yap) \neq 0$, the $y$-coordinate along $\gamap$
firstly monotonically increases from $y_0$ to $\yap$, reaches the maximum $y=\yap$,
and then monotonically decreases from $\yap$ to $y_0$.

2. If  $\yap < \om^+$ and $a'(\yap) = 0$, the $y$-coordinate on $\gamap$ always monotonically increases,
and $\gamap$ asymptotically tends to the horizontal line $y=\yap$, which is also a geodesic.

3. If  $\yap = \om^+$, the $y$-coordinate on $\gamap$ monotonically increases from $y_0$ to $\om^+$,
and the geodesic $\gamap$ leaves the domain $\Om$ through the boundary $y = \om^{\pm}$.

Similar statements hold true for $\gamam \in \Gamo$.
\end{theorem}

{\bf Proof}.
For definiteness, consider $\gamap$ (for $\gamam$ the proof is similar). A necessary condition for $\gamap$ to return to the initial line $y=y_0$ is that $\gamap$ contains a point $(x,y)$, $y_0 < y < \om^+$, with the horizontal tangential direction $p=0$. The case $\ha^2 = 0, \infty$ is trivial: substituting $p=0$ in the equation \eqref{19} with $h^2 = 0, \infty$, we obtain $a(y)=0$. This equation has no solutions in $\Om$, hence $\gamap$ with $\ha^2 = 0, \infty$ do not return to the initial line $y=y_0$.

Now assume that $\ha^2 \neq 0, \infty$.
Substituting $p=0$ in \eqref{19}, we get the equation $a(y)(a(y)-\ha^2)=0$, which is equivalent to $a(y)=\ha^2$, since $a(y)$ does not vanish in $\Om$. Therefore $\yap < \om^+$ is a necessary condition for $\gamap$ to return to the initial line $y=y_0$.

\begin{figure}[h]
\begin{center}
\includegraphics{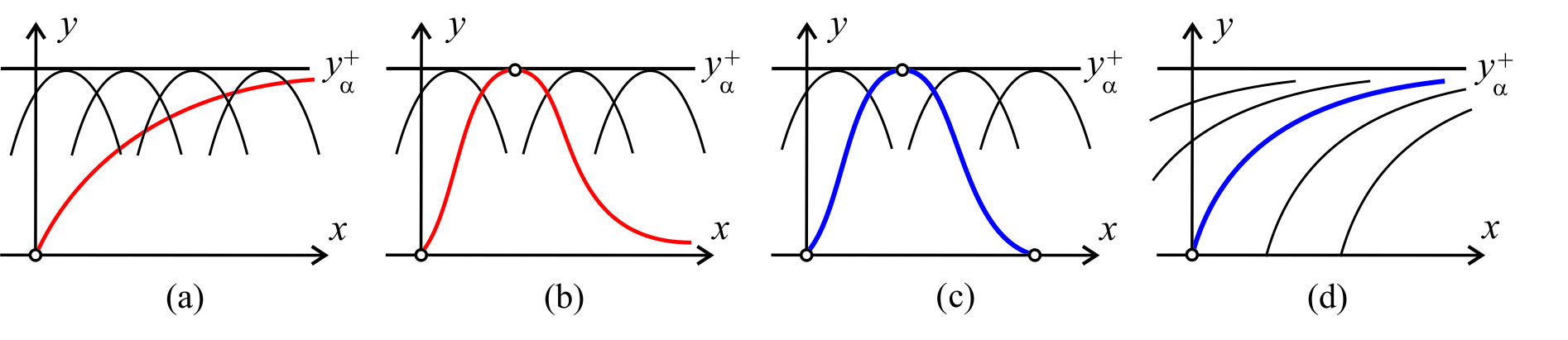}
\caption{
The case~1 in Theorem~\ref{T6}: impossible behavior of $\gamap$ (red line) on (a) and (b);
the only possible situation ($\gamap$ is the blue line) on (c).
The case~2 in Theorem~\ref{T6}: $\gamap$ is the blue line presented on (d).
Here the point $q_0$ coincides with the origin.
}\label{Fig2}
\end{center}
\end{figure}

{\bf The case 1}: $\yap < \om^+$ and $a'(\yap) \neq 0$.
By lemma~\ref{L4}, in a neighborhood of any point of the line $y=\yap$, the equation \eqref{19} with $h^2 = \ha^2$ has the family of non-singular solutions
\begin{equation}
y=\yap+f(x-c), \ \  f(0)=f'(0)=0, \ \, f''(0)<0,  \quad c \in \bR,
\label{24}
\end{equation}
which are geodesics, and the singular solution $y=\yap$, which is not a geodesic.
(Non-singular solutions \eqref{24} correspond to parabolas ${\ti y}= \frac{1}{4} ({\ti x}-c)^2$ in the normal form
${\ti p}^2 = {\ti y}$, and $y=\yap$ corresponds to ${\ti y}=0$.)
By lemma~\ref{L5}, $\gamap$ tends arbitrarily close to the line $y=\yap$, hence
it coincides with one of non-singular solutions \eqref{24}, and reaches the line $y=\yap$ at a finite point $(x_*,\yap)$. After passing through the point $(x_*,\yap)$, the $y$-coordinate along $\gamap$ monotonically decreases from $\yap$ to $y_0$. By Assumption~\ref{AS2}, $\gamap$ reaches the horizontal line $y=y_0$ at a finite point.
It is presented on Fig.~\ref{Fig2} (c), while on Fig.~\ref{Fig2} (a,b) two impossible situations are presented.

{\bf The case 2}: $\yap < \om^+$ and $a'(\yap) = 0$.
Since $a(\yap)=\ha^2$, on the line $y=\yap$ the equation \eqref{19} with $h^2 = \ha^2$ is equivalent to $p^2=0$,
and $\gamap$ cannot intersect the line $y=\yap$ transversally.
On the other hand, from $a'(\yap) = 0$ it follows that the line $y=\yap$ is a geodesic (Lemma~\ref{L4}),
and $\gamap$ cannot be tangent to the line $y=\yap$ at a finite point (if this case there are two geodesics passing through the same point with the same tangential direction).
Thus there is only one possible situation: $\gamap$ asymptotically tends to the horizontal geodesic $y=\yap$ and
it never returns to the initial line $y=y_0$, as it is presented on Fig.~\ref{Fig2} (d).
\hfill\ebox

\begin{remark}
\label{R4}
For $\ha^2 = 0, \infty$, we have $\yapm = \om^{\pm}$. Theorem~\ref{T6} asserts that geodesics $\gamapm \in \Gamo$ with
$\ha^2 = 0, \infty$ (if they exist) never return to the initial line $y=y_0$.
If the metric \eqref{18} is diagonal, i.e., $b(y) \equiv 0$, geodesics $\gamapm$ with $\ha^2=0$ are the halves of vertical lines $x = \const$, while the $x$-coordinate along each $\gamap$, $\ha^2 \neq 0$, monotonically increases or decreases.
\end{remark}

\begin{figure}[h]
\begin{center}
\includegraphics{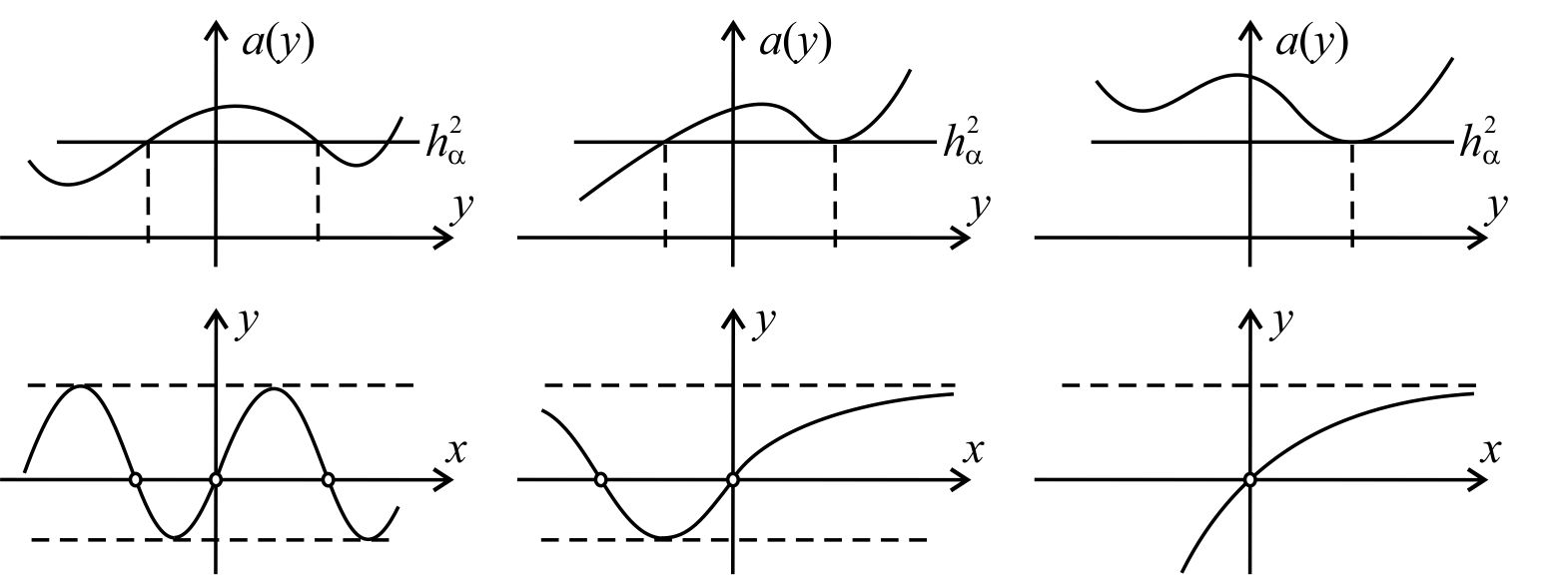}
\caption{Examples of functions $a(y)$ (up) and geodesics $\gamapm \in \Gamo$, $q_0=0$ (down).
}\label{Fig3}
\end{center}
\end{figure}


\subsection{Riemannian and Lorentzian metrics}\label{R-and-L}

For illustration, consider the cases when the metric \eqref{18} is Riemannian or Lorentzian on the whole $(x,y)$-plane. Then geodesics $\gamapm \in \Gamo$ can be indexed by their tangential directions
$p \in \obR = \bR \cup \infty$ at $q_0$, and $\ha^2$ can be expressed through $p$ after substituting $y=y_0$ in \eqref{19} or \eqref{20}. However, in order to obtain agreement with previous notation used for pseudo-Riemannian metrics, we shall put $\alf$ equal to the value $p^{-1}$ at $q_0$.

\begin{figure}[h!]
\begin{center}
\includegraphics{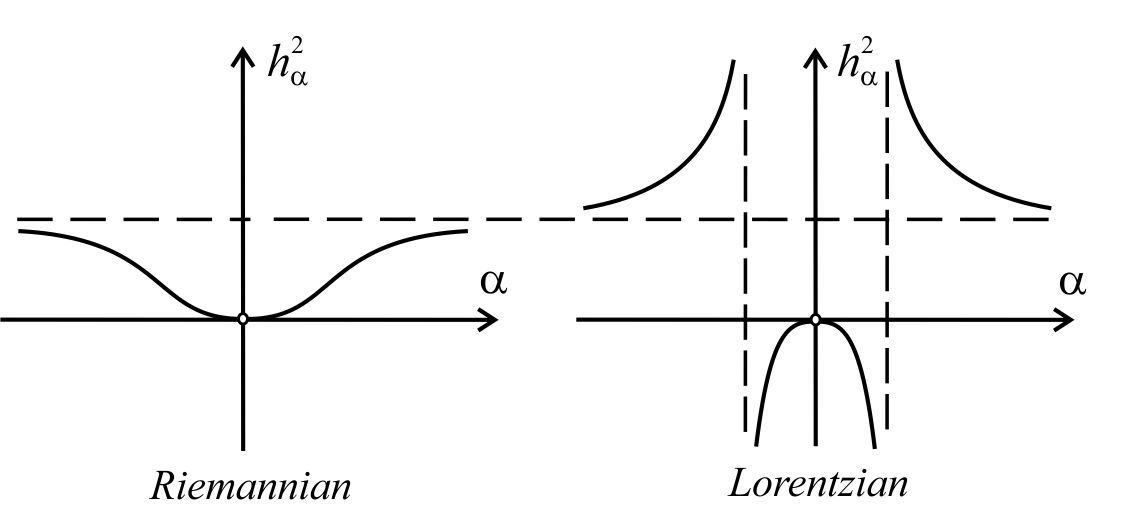}
\caption{The dependence $\ha^2$ on $\alf$ for a Riemannian metric (left) and a Lorentzian metric (right). }
\label{Fig4}
\end{center}
\end{figure}

Without loss of generality we shall further assume that $a(y)>0$ for all $y$ and $b(0) = 0$
(this can be obtained by the change of variables $x \to x + b(0)/a(0) y$).
The formula \eqref{20} implies that the set $\HH$ of all possible values of the constant $h^2$ is defined by the formula
\begin{equation}
\HH = \biggl\{ \ha^2 = \frac{(\alf a_0)^2}{\alf^2 a_0 + c_0}, \quad  \alf\in \obR = \bR \cup  \infty \biggr\},
\quad \textrm{where}  \ \, a_0=a(0), \ \, c_0=c(0).
\label{25}
\end{equation}
Formula \eqref{25} shows that $\HH = [0,a_0]$ for Riemannian metrics ($c_0>0$) and $\HH = [a_0,+\infty) \cup (-\infty,0] \cup \infty$ for Lorentzian metrics ($c_0<0$); see Fig.~\ref{Fig4}.

\subsubsection{Riemannian metrics}\label{Riemannian}

By Theorem \ref{T6}, to define the behavior of geodesics $\gamapm \in \Gamo$ with $\ha^2 \neq 0,\infty$,
it is sufficient to consider the equation $a(y) = \ha^2$ with all possible values $0 < \ha^2 \le a_0$.
This gives the following conclusions:

\begin{itemize}
\item
Geodesics $\gam^{\pm}_{\alf}$ and $\gam^{\pm}_{-\alf}$ with the same superscript return (or do not return)
to the initial horizontal line $y=y_0$ simultaneously.

\item
If $a(y) \ge a_0$ for all $y$, then all geodesics $\gamapm \in \Gamo$ with $\ha^2 \neq a_0$ (i.e., $\alf \neq \infty$) do not return to $y=y_0$, while $\ha^2 = a_0$ correspond to the geodesic $\gam_{\infty}$ being the the horizontal line $y = y_0$. Otherwise the family $\Gamo$ contains an infinite number of geodesics returning to $y=y_0$.

\item
If $a'(y_0) \neq 0$, then the family $\Gamo$ contains an infinite number of geodesics returning to the line $y=y_0$. Moreover, if $a'(y_0) < 0$ $(>0)$, there exists $\eps>0$ such that all geodesics $\gamap \in \Gamo$ (respectively, $\gamam \in \Gamo$) with $0< |\alf|^{-1} < \eps$ return to the line $y=y_0$.

\item
If $a(y)$ has a strict local maximum at $y_0$, then the horizontal line $y = y_0$ is a geodesic (with the constant $h_{\infty}^2=a_0$), and there exists $A>0$ such that all $\gamapm \in \Gamo$ with $A < |\alf| < \infty$ return to the line $y=y_0$.
\end{itemize}

\subsubsection{Lorentzian metrics}

For Lorentzian metrics \eqref{18} with $a(y)>0$, let us separate the family $\Gamo$ into three subfamilies: $\Gamo^t$, $\Gamo^s$, $\Gamo^i$ containing timelike, spacelike, and isotropic geodesics, respectively.
Accordingly, the set $\HH$ also consists of three parts:
$$
\HH = [a_0,+\infty) \cup (-\infty,0] \cup \infty,
$$
which correspond to $\Gamo^t$, $\Gamo^s$, $\Gamo^i$, respectively.
Theorem \ref{T6} asserts that geodesics $\gamapm \in \Gamo^s \cup \Gamo^i$ do not return to the line $y=y_0$,
and the $y$-coordinate monotonically increases from $y_0$ to $+\infty$ on $\gamap$ and monotonically decreases from $y_0$ to $-\infty$ on $\gamam$.

Geodesics $\gamapm \in \Gamo^t$ can be studied similarly to geodesics in Riemannian metrics,
with evident changes. It is sufficient to consider the equation $a(y) = \ha^2$ with all possible values
$\ha^2 \ge a_0$. The properties established for Riemannian metrics above, can be transferred to Lorentzian metrics with the following alterations: $\Gamo$ should be replaced by $\Gamo^t$, maximum should be replaced by minimum, and
inequality symbols in $a(y) \ge a_0$, $a'(y_0) > 0$, $a'(y_0) < 0$ should be reversed.

\subsection{Pseudo-Riemannian metrics}\label{Pseudo-Riemannian}

Consider smooth metric \eqref{18} that is Riemannian in the lower semiplane $y<y_0$ and Lorentzian in the upper semiplane $y>y_0$. The horizontal line $y=y_0$ is the degeneracy curve $\AA$.
Without loss of generality assume that $q_0=0$ (the origin), and $a(0)=a_0>0$, $b(0)=c(0)=0$.
Also assume that the following genericity conditions hold:
$$
a_1:= a'(0) \neq 0,  \quad  c_1:= \frac{4}{9} c'(0) \neq 0.
$$
The second condition gives $c_1 < 0$, since the upper semiplane is Lorentzian.

If $a_1 > 0$, the isotropic direction $p = \infty$ is the only admissible direction at $0$. If $a_1 < 0$, it is supplemented by two timelike admissible directions $p_{1,2} = \pm \frac{2}{3}\sqrt{a_1/c_1}$.
Geodesics with admissible directions $p_{1,2}$ can be studied similarly to what was done for geodesics in Riemannian metrics (substituting $y=0$, $b=c=0$, $p=p_{1,2}$ in \eqref{19} or \eqref{20}, one can see that the corresponding constant $h^2 = a_0$). Whence we focus our attention on the family $\Gamo = \{ \gamapm \}$ of geodesics outgoing from $0$ with the isotropic tangential direction $p = \infty$.

As in the case of Lorenzian metrics, we represent the family $\Gamo$ as the union $\Gamo = \Gamo^t \cup \Gamo^i \cup \Gamo^s$. Moreover, we separate the subfamily $\Gamo^t$ into two parts: $\Gamo^{t,L}$ and $\Gamo^{t,R}$,
which contain timelike geodesics $\gamap$ lying in the Lorentzian semiplane and
timelike geodesics $\gamam$ lying in the Riemannian semiplane, respectively.
To establish the ranges of $\ha^2$ for each family of $\Gamo$, one can use the local coordinates from Theorem~\ref{T4}. Then from \eqref{10} and \eqref{20} we obtain
$$
\ha^2 =
\lim_{\tau \to 0} \frac{(a(y) + b(y)p)^2}{a(y) + 2b(y)p + c(y)p^2} =
\lim_{\tau \to 0} \frac{(a_0+O(\tau))^2}{a_0 + \alf^{-2}c_1 Y_{\alf}(0) +O(\tau)} =
\frac{(\alf a_0)^2}{\alf^{2}a_0 + c_1 Y_{\alf}(0)},
$$
where $Y_{\alf}(0) = \pm 1$ for $\gamapm$, respectively.
This yields the results presented in Table~1.

\begin{table}[h]
\label{Table1}
\begin{center}
\begin{tabular}{|c|c|c|c|c|c|}
\hline
Geodesic &  Type & Semiplane & Range of $\alf$ & Range of $\ha^2$   \\
\hline
$\gamam \in \Gamo^{t,R}$ & timelike & $y<0$ & $\alf \in \bR$     & $0 \le \ha^2 < a_0$   $\phantom{\Bigl|}$  \\
\hline
$\gamap \in \Gamo^{t,L}$ & timelike & $y>0$ & $|\alf| > \sqrt{|c_1|/a_0}$ & $ a_0 < \ha^2 < +\infty$ $\phantom{\Bigl|}$ \\
\hline
$\gamap \in \Gamo^{i}$   & isotropic & $y>0$ & $\alf = \pm \sqrt{|c_1|/a_0}$ & $\ha^2 = \infty$ $\phantom{\Bigl|}$ \\
\hline
$\gamap \in \Gamo^{s}$   & spacelike & $y>0$ & $|\alf| < \sqrt{|c_1|/a_0}$ & $-\infty < \ha^2 \le 0$  $\phantom{\Bigl|}$ \\
\hline
\end{tabular}
\caption{Geodesics outgoing from a generic parabolic point.}
\end{center}
\end{table}

\medskip

Theorem \ref{T6} asserts that spacelike and isotropic geodesics, and likewise, the timelike geodesic $\gam_0^-$ do not return to the initial line $y=0$, since the $y$-coordinate monotonically increases from $0$ to $+\infty$ on $\gamap$ or decreases from $0$ to $-\infty$ on $\gamam$.
Remaining timelike geodesics $\gamam \in \Gamo^{t,R}$, $\alf \neq 0$, and $\gamap \in \Gamo^{t,L}$
can be studied similarly to geodesics in Riemannian and Lorentzian metrics, respectively.
Namely, for $\gamam \in \Gamo^{t,R}$ we need to consider the equation $a(y) = \ha^2$
on the interval $y<0$ with all possible constants $0<\ha^2<a_0$.
For $\gamap \in \Gamo^{t,L}$ we need to consider the equation $a(y) = \ha^2$
on the interval $y>0$ with all possible constants $a_0<\ha^2<+\infty$.

Observe some corollaries:

\begin{itemize}
\item
Geodesics $\gamma^{+}_{\pm \alf}$ $(\gamma^{-}_{\pm \alf})$ with the same superscript
(being branches of the same semicubic parabola)
return (or do not return) to the initial line $y=0$ simultaneously.

\item
If $a'(0)>0$, then for any $\eps>0$ there exists $A>0$ such that all timelike geodesics
$\gamamp \in \Gamo^{t}$ with $|\alf|>A$ do not leave the strip $|y| < \eps$ and return to the line $y=0$.

\item
If $a'(y)>0$ for all $y$ and $\inf a(y) = 0$, $\sup a(y) = +\infty$, then all timelike geodesics $\gamamp \in \Gamo^{t}$ with $\alf \neq 0$ return to the line $y=0$.

\item
If $a'(0)<0$, then there exists $\eps>0$ such that all geodesics leave the strip $|y| < \eps$.

\item
If $a'(y) \le 0$ for all $y$, then all geodesics do not return to the line $y=0$.
\end{itemize}

\medskip

\begin{figure}[h]
\begin{center}
\includegraphics{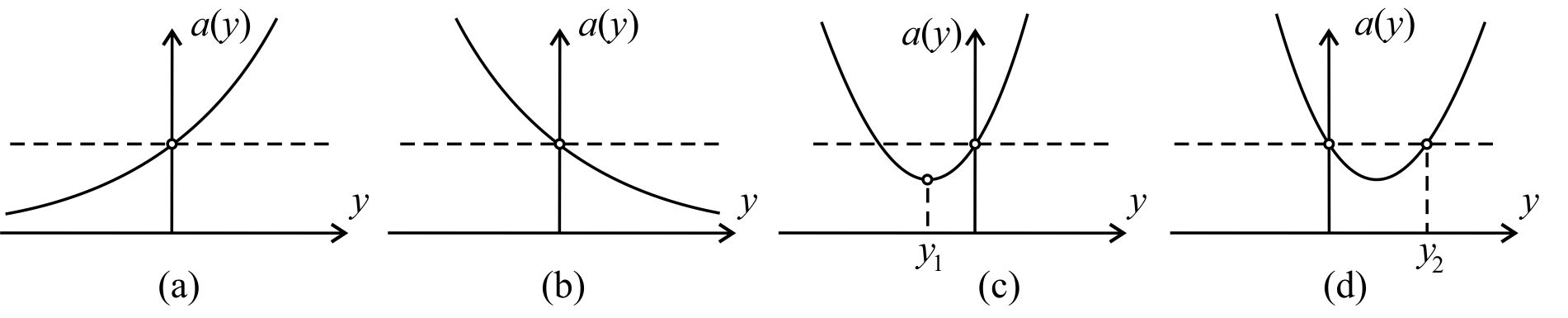}
\caption{Four examples of $a(y)$. The correspond families $\Gamo$ are presented on Fig.~\ref{Fig5},~\ref{Fig14}.
}\label{Fig13}
\end{center}
\end{figure}

{\bf Example 3.2.}
Consider the functions $a(y)$ presented on Fig.~\ref{Fig13}. In the cases (a) and (b) the family $\Gamo$ is presented on Fig.~\ref{Fig5}, left and right, respectively.
In the case (a) all geodesics $\gamapm \in \Gamo^t$, $\ha^2 \neq 0$, return to the line $y=0$, while
in the case (b) the family $\Gamo$ does not contain geodesics returning to the line $y=0$ at all.

\begin{figure}[h]
\begin{center}
\includegraphics{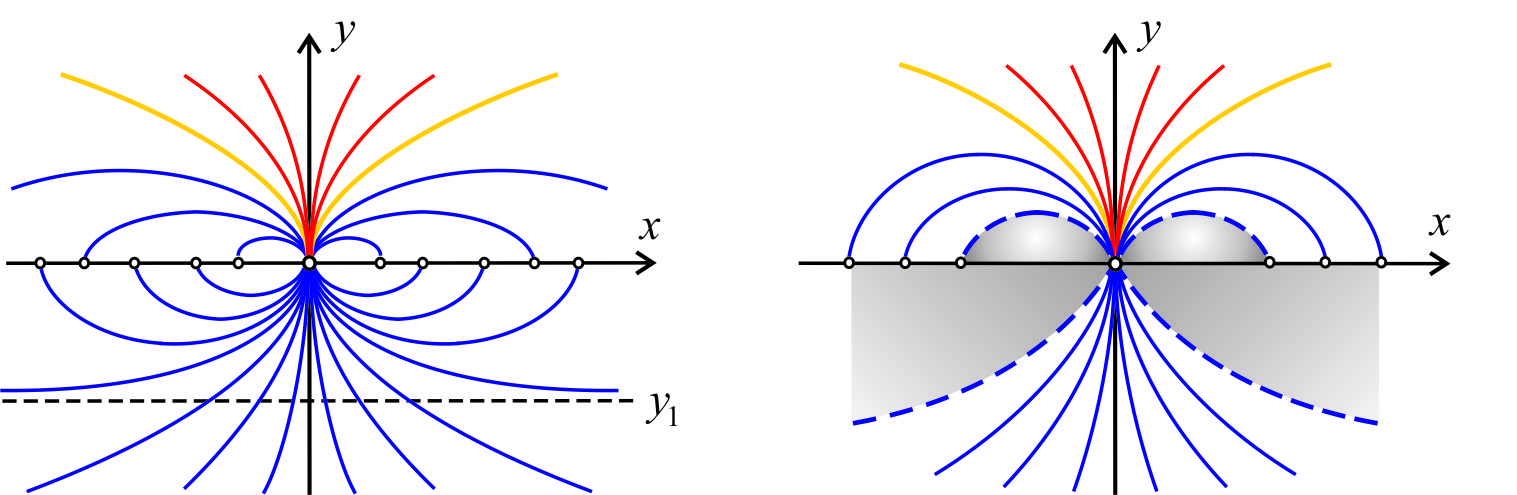}
\caption{The family $\Gamo$ for the functions $a(y)$ presented on Fig.~\ref{Fig13} (c), (d).
Timelike, spacelike, and isotropic geodesics are depicted as the blue, red, and yellow lines respectively.
The horizontal geodesic $y=y_1$ that does not pass through $0$ (left) and geodesics passing through $0$ with non-isotropic admissible directions (right) are depicted as the dashed lines.
}\label{Fig14}
\end{center}
\end{figure}

In the case (c) the family $\Gamo$ is presented on Fig.~\ref{Fig14} (left).
For every $\ha^2 > a(0)$ we have $\yap < +\infty$ and $a'(\yap) \neq 0$,
hence all geodesics $\gamap \in \Gamo^{t,L}$ return to the line $y=0$.
The function $a(y)$ has a global minimum at $y_1<0$.
Hence the horizontal line $y=y_1$ is a geodesic with $\ha^2 = a(y_1)$, and
$\gamam \in \Gamo^{t,R}$ have different behavior:
geodesics $\gamam$ with $a(y_1) < \ha^2 < a(0)$ return to $y=0$,
geodesics $\gamam$ with $0 < \ha^2 < a(y_1)$ do not return to $y=0$ (the $y$-coordinate decreases from $0$ to $-\infty$), and the sole geodesic $\gamam$ with $\ha^2 = a(y_1)$ asymptotically tends to $y=y_1$.

In the case (d) the family $\Gamo$ is presented on Fig.~\ref{Fig14} (right).
Here there exist two supplemented geodesics $\delta^{\pm}_{1,2} \notin \Gamo$
passing through $0$ with non-isotropic admissible directions $p_{1,2} = \pm \frac{2}{3}\sqrt{a_1/c_1}$.
The corresponding constant $h^2 = a(0)$.
The geodesics $\delta^{\pm}_{1,2}$ are depicted as the dashed lines, the domains between them and the $x$-axis (colored in grey) do not contain any geodesics from $\Gamo$. For every $\ha^2 > a(0)$ we have $\yap < +\infty$ and $a'(\yap) \neq 0$, hence all geodesics $\gamap \in \Gamo^{t,L}$ return to the line $y=0$.
For every $0 < \ha^2 < a(0)$ the value $\yam = -\infty$,
hence all geodesics $\gamam \in \Gamo^{t,R}$ do not return to the line $y=0$.

\hfill\ebox


\subsection{Two geometric examples: surfaces of revolution}\label{Examples}

An example of the pseudo-Riemannian metric \eqref{18} is a metric induced on a surface of revolution $\M$
embedded in three-dimensional Minkowski space $\bR^3_1$ with signature $(++-)$.
Here by surface of revolution we mean a surface that is invariant with respect to Lorenzian transformations of the space $\bR^3_1$. In Cartesian coordinates $r=(r_1,r_2,r_3)$, the metric in $\bR^3_1$ has the form $ds^2 = dr_1^2+dr_2^2-dr_3^2$,

\subsubsection{Geodesics on a sphere in Minkowski space}

Let $\M$ be a (Euclidean) sphere in $\bR^3_1$ given by the equation $r_1^2 + r_2^2 + r_3^2 = 1$.
The metric induced on $\M$ has transverse parabolic points, which form two circles (parallels) $C_N$ and $C_S$ being the intersections of the sphere with the planes $r_3 = \pm \frac{1}{\sqrt{2}}$ (the subscripts $N=North$ and $S=South$ correspond to the plus and the minus, respectively).
The parallels $C_N$, $C_S$ separate the sphere $\M$ into three regions, where the metric has constant signatures.
The North region $\M_N \colon r_3> \frac{1}{\sqrt{2}}$ and the South region $\M_S \colon r_3< -\frac{1}{\sqrt{2}}$ are
Riemannian, while the equatorial region $\M_E \colon |r_3|< \frac{1}{\sqrt{2}}$ is Lorenzian.

To bring the metric on the sphere to the form \eqref{18}, one can use the spherical coordinates
$r_1 = \sin \th \cos \phi$, $r_2 = \sin \th \sin \phi$, $r_3 = \cos \th$,
$0 \le \th \le \pi$, $\phi \in \bR$.
Then $ds^2 = \sin^2 \th \, d\phi^2 + \cos{2\th} \,d\th^2$.
Moreover, it is convenient to use new coordinated $x = \sqrt{2} \phi$, $y = 2\th - \frac{\pi}{2}$.
Multiplying the metric by a constant factor, we have
\begin{equation}
ds^2 = (1+\sin y)\, dx^2 - \sin y\, dy^2.
\label{26}
\end{equation}
Then the equator $E$ is given by $y=\frac{\pi}{2}$, and the parallels $C_N$, $C_S$ are given by $y = 0$, $y = \pi$, respectively. Since the Gaussian curvature of sphere is everywhere positive, at any parabolic point
there exists a unique admissible (isotropic) direction $p = \infty$.

Remark that the metric \eqref{26} is degenerate not only on $C_N$ and $C_S$,
but also at the North and South poles: $y=-\frac{\pi}{2}, \frac{3\pi}{2}$.
However, singularities at the poles appear only due to the singularity of the spherical coordinate system and do not have a geometric meaning. In a neighborhood of each pole, there exist local coordinates where the metric is Riemannian, and at each pole the family $\Gamo$ consists of meridians on the sphere.

\medskip

\begin{figure}[h!]
\begin{center}
\includegraphics{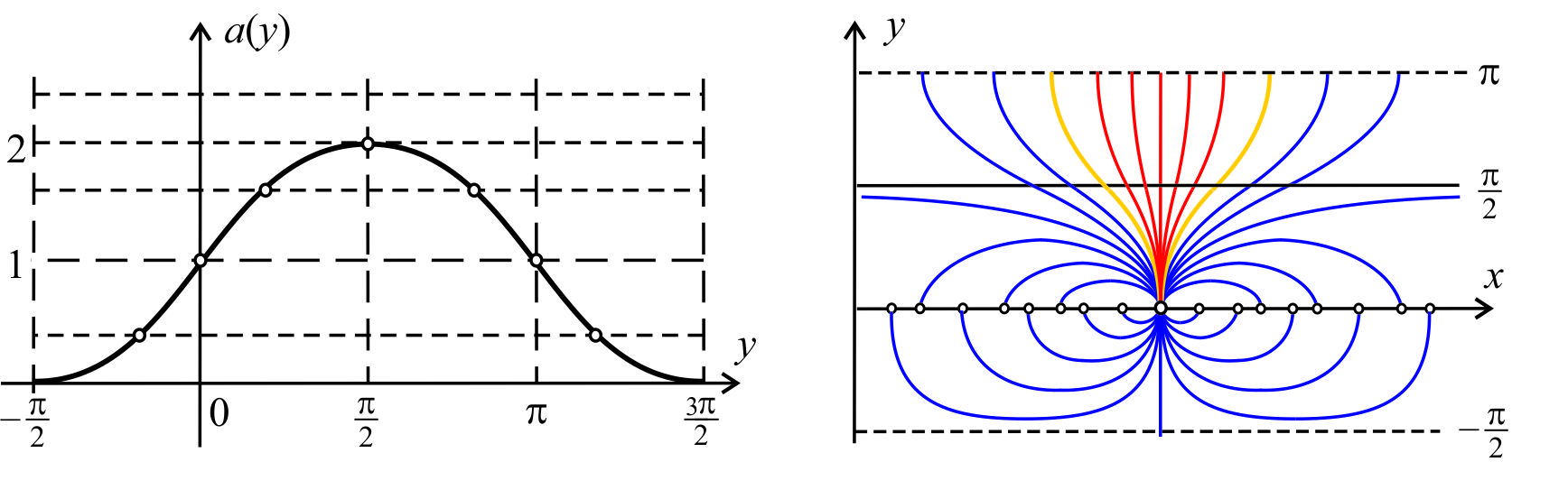}
\caption{The graph of the function $a(y)=1+\sin{y}$ (left).
Geodesics outgoing from a parabolic point $q_0 \in C_N$ on the $(x,y)$-plane (right).
Timelike, spacelike, and isotropic geodesics are depicted as the blue, red, and yellow lines respectively.
}\label{Fig6}
\end{center}
\end{figure}

\begin{figure}[h!]
\begin{center}
\includegraphics{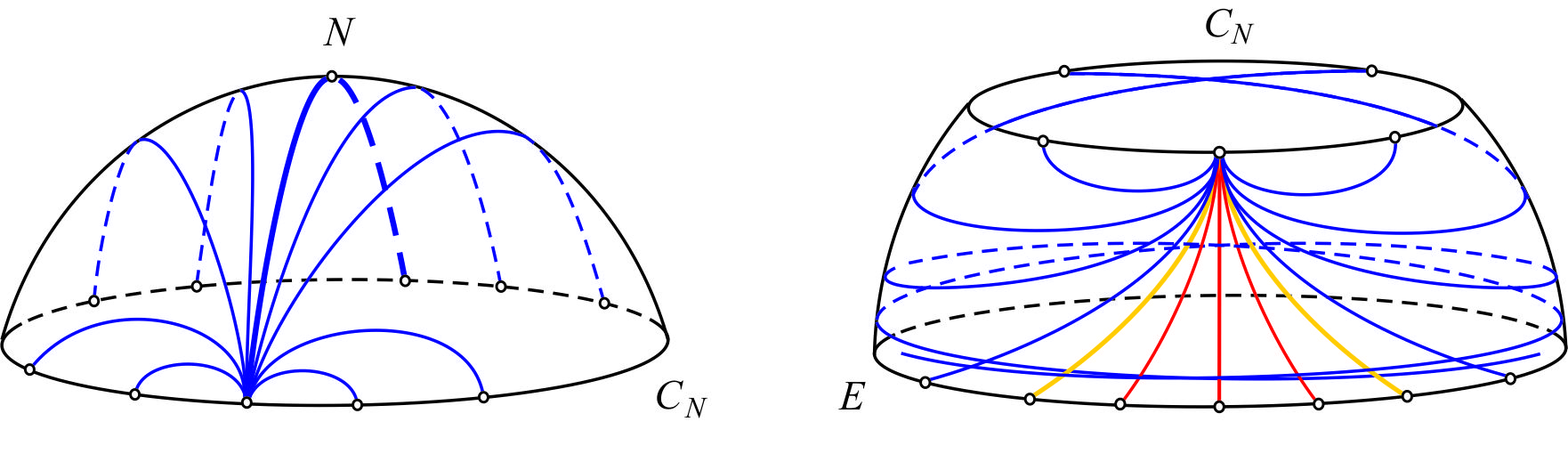}
\caption{On the left: the region $\M_N$ and geodesics $\gamam \in \Gamo$ for $q_0 \in C_N$
On the right: the Northern half of the region $\M_E$ and geodesics $\gamap \in \Gamo$ for $q_0 \in C_N$.
Timelike, spacelike, and isotropic geodesics are depicted as the blue, red, and yellow lines respectively.
}
\label{Fig7}
\end{center}
\end{figure}

To analyze geodesics on the whole sphere, it is sufficient to consider the families $\Gamo$ in turn for points
$q_0 \in \M_N, C_N, \M_E$. The corresponding domain $\ov \Om$ is $\M_N$ ($-\frac{\pi}{2}<y<0$), $\M_N \cup \M_E$ ($-\frac{\pi}{2}<y<\pi$), $\M_E$ ($0<y<\pi$), respectively. Applying Theorem~\ref{T6}, we get the following results (see Fig.~\ref{Fig6}, \ref{Fig7}):

\begin{itemize}
\item
The equator is a unique geodesics (with $\ha^2=2$) that does not approach the set of parabolic points $C_N \cup C_S$ arbitrarily closely. Consequently, it is a unique geodesics whose natural parametrization does not have singularity.

\item
All regions $\M_N, \M_E, \M_S$ contain geodesics with $\ha^2=0$ -- bows of meridians on the sphere,
which are timelike in $\M_N$, $\M_S$ and spacelike in $\M_E$.

\item
In the Riemannian region $\M_N$ ($\M_S$), all remaining geodesics ($0< \ha^2 < 1$) are curves of finite length with two cusps at endpoints on the parallel $C_N$ (respectively, $C_S$).

\item
In the Lorenzian region $\M_E$, all remaining geodesics are curves of finite length with two cusps at endpoints ($\ha^2 \neq 2$) or curves of infinite length with one endpoint ($\ha^2=2$) winding round the equator.
More exactly, there are seven different classes of geodesics in $\M_E$ enumerated in Table~2.
\end{itemize}

\begin{table}[h]
\label{Table2}
\begin{center}
\begin{tabular}{|c|c|c|c|c|c|}
\hline
{} & Type & Range of $\ha^2$ & 1st endpoint & 2nd endpoint &   Brief description  $\phantom{\Bigl|}$  \\
\hline
1  & timelike &  $1 < \ha^2 < 2$  & cusp on $C_N$ & cusp on $C_N$  & Belong to the $N$-hemisphere $\phantom{\Bigl|}$  \\
\hline
2  & timelike &  $\ha^2 = 2$  & cusp on $C_N$ &  ---  & Start from $C_N$ and do not return,    $\phantom{\Bigl|}$ \\
  &  &   &  &                                          & wind round $E$ in the $N$-hemisphere \\
\hline
3  & timelike &  $2 < \ha^2 < +\infty$  & cusp on $C_N$ & cusp on $C_S$  & Intersect $E$, connect $C_N$ and $C_S$    $\phantom{\Bigl|}$ \\
\hline
4  & isotropic &  $\ha^2 = \infty$  & cusp on $C_N$ & cusp on $C_S$  & Intersect $E$, connect $C_N$ and $C_S$   $\phantom{\Bigl|}$ \\
\hline
5  & spacelike &  $-\infty < \ha^2 < 0$  & cusp on $C_N$ & cusp on $C_S$  & Intersect $E$, connect $C_N$ and $C_S$   $\phantom{\Bigl|}$ \\
\hline
6  & timelike &  $1 < \ha^2 < 2$  & cusp on $C_S$ & cusp on $C_S$  & Belong to the $S$-hemisphere   $\phantom{\Bigl|}$ \\
\hline
7  & timelike &  $\ha^2 = 2$  & cusp on $C_S$ &  ---  & Start from $C_S$ and do not return,    $\phantom{\Bigl|}$ \\
  &  &   &  &                                          & wind round $E$ in the $S$-hemisphere \\
\hline
\end{tabular}
\caption{Seven classes of geodesics (with the exception of the equator and meridians) in $\M_E$.}
\end{center}
\end{table}

We present the detailed analysis of the family $\Gamo$ for $q_0 \in C_N$ only.
$\Gamo$ contains two geodesics $\gam^+_0 \in \Gamo^s$ and $\gam^-_0 \in \Gamo^{t,R}$
with $\ha^2=0$, which are the bows of the meridian passing through the point $q_0$.
Further we analyze geodesics $\gamapm \in \Gamo$ with $\ha^2 \neq 0$.

For $\gamam \in \Gamo^{t,R}$, consider the equation $1+\sin{y} = \ha^2$
on the interval $-\frac{\pi}{2}< y < 0$ with all possible constants $0<\ha^2<1$.
Given $\ha^2$, it has a unique solution $\yam$, and the condition $a'(\yam) \neq 0$ holds (Fig.~\ref{Fig6}, left).
Hence each geodesic $\gamam \in \Gamo^{t,R}$ goes out from the point $q_0$ toward the North,
turns back on the parallel $y=\yam$, and returns to the initial parallel $C_N$.
Geodesics of the subfamily $\Gamo^{t,R}$ fill the region $\M_N$ as it is presented on Fig.~\ref{Fig7} (left).

Now consider geodesics $\gamap \in \Gamo^{t,L}$.
We have to distinguish three different cases: $1 < \ha^2 < 2$, $\ha^2 = 2$, and $2 < \ha^2 < +\infty$,
which correspond to the classes 1, 2, and 3 from Table~2.

For $1 < \ha^2 < 2$, we have $0< \yap < \frac{\pi}{2}$ and $a'(\yap) \neq 0$.
Hence each geodesic $\gamap \in \Gamo^{t,L}$ with $1 < \ha^2 < 2$ goes out from $q_0$ toward the South,
turns back on the parallel $y=\yap$, and returns to the initial parallel $C_N$.

For $\ha^2 = 2$, we have $\yap = \frac{\pi}{2}$ and $a'(\yap) = 0$.
Both geodesics $\gamap \in \Gamo^{t,L}$ with $\ha^2 = 2$ go out from $q_0$ toward the South,
but in contrast to the previous case, they do not return back to the parallel $C_N$.
They wind round the equator staying in the North hemisphere.

For $2 < \ha^2 < +\infty$, we have $\yap = +\infty$. Hence all geodesics $\gamap \in \Gamo^{t,L}$
with $2 < \ha^2 < +\infty$ do not return back.
The  $y$-coordinate on each of them monotonically increases from $0$ to $\pi$, and each
geodesic intersects the equator and reaches the parallel $C_S$ at a certain point.

Finally, consider geodesics $\gamap \in \Gamo^{i}$ and $\gamap \in \Gamo^{s}$
(the classes 4, 5 from the Table~2.)
As in the previous case, $\yap = +\infty$, and the behavior of geodesics
$\gamap \in \Gamo^{i}$ and $\gamap \in \Gamo^{s}$
is similar to $\gamap \in \Gamo^{t,L}$ with $2 < \ha^2 < +\infty$.
Geodesics of all these classes are presented on Fig.~\ref{Fig6} (right) and Fig.~\ref{Fig7} (right).


\subsubsection{Geodesics on a torus in Minkowski space}

Let $\M$ be a torus, whose axis of revolution coincides with the $r_3$-axis of the ambient Minkowski space.
The group of symmetries of the metric induced on $\M$ includes Euclidean rotations of the $(r_1,r_2)$-plane.

Without loss of generality assume that the torus is given by the formula
$$
r_1 = (\rho + \cos y) \cos x, \ \ r_2 = (\rho + \cos y) \sin x, \ \
r_3 = \sin y,
$$
where $x, y \in \bR$ (mod $2\pi$) are the coordinates on $\M$, and $\rho > 1$ is a constant determining size of
the torus. By $N$ and $S$ denote the North and South parallels of the torus: $y=\pm \frac{\pi}{2}$.
They separate the torus into the outer part $\M^+ \colon |y|<\frac{\pi}{2}$
and the inner part $\M^- \colon |y-\pi|<\frac{\pi}{2}$,
with positive and negative Gaussian curvature, respectively.
The metric induced on $\M$ reads
\begin{equation}
ds^2 = (\rho+\cos y)^2 \, dx^2 - \cos{2y}\, dy^2.
\label{27}
\end{equation}

The metric \eqref{27} has transverse parabolic points, which form four parallels $C_N^+$, $C_N^-$, $C_S^-$, $C_S^+$, given by $y=\frac{\pi n}{4}$, $n=1,3,5,-1$, respectively.\footnote{
The subscript $N$ ($S$) indicates that the given parallel belongs to the North (South) half of the torus,
while the superscript $\pm$ indicates that it belongs to the outer or inner part of the torus, in accordance
with the sign of the Gaussian curvature.
}
A parabolic point $q_0 \in C_N^+ \cup C_S^+$ has a unique admissible (isotropic) direction $p=\infty$,
while $q_0 \in C_N^- \cup C_S^-$ has three admissible directions:
$p=\infty$ and two timelike directions $p_{1,2} = \pm \frac{\sqrt[4]{2}}{2}$,
to which correspond geodesics with $h^2 = (\rho-\frac{1}{\sqrt{2}})^2$.

\begin{figure}[h]
\begin{center}
\includegraphics{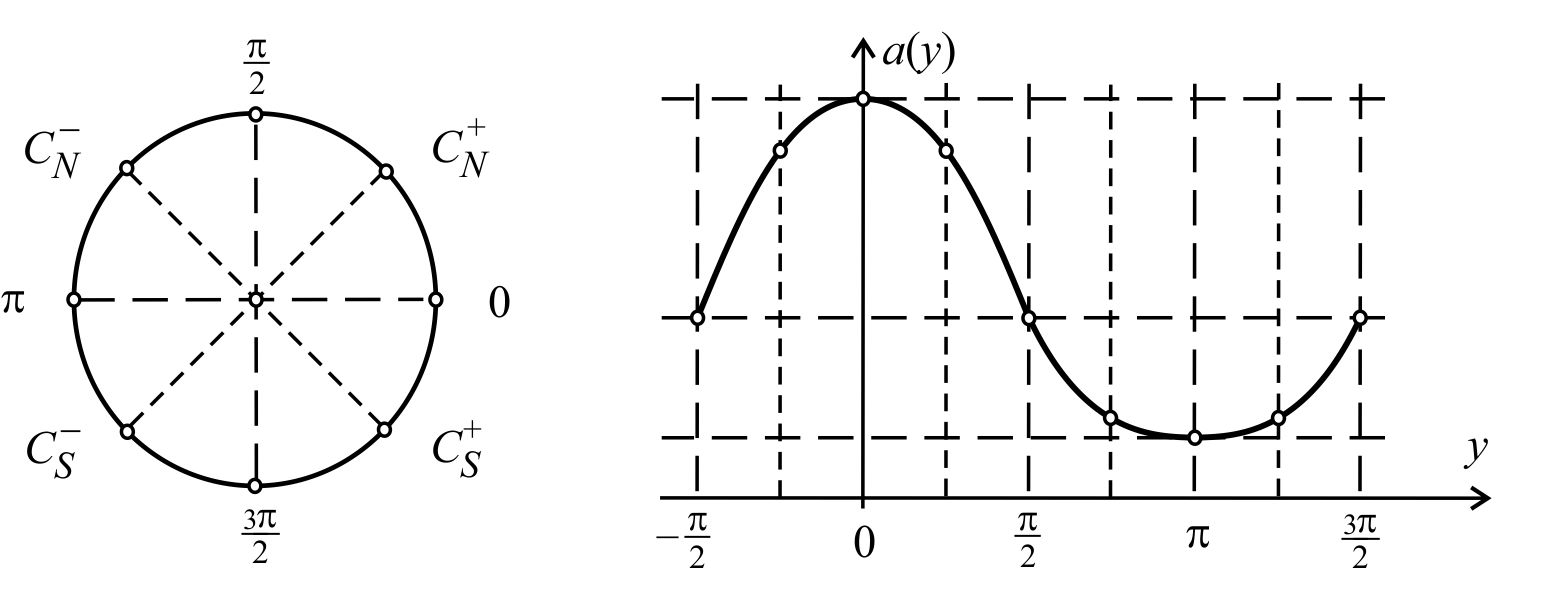}
\caption{
On the left: section of the torus along the meridian $x=0$.
On the right: the graph of the function $a(y)=(\rho+\cos y)^2$.
}
\label{Fig8}
\end{center}
\end{figure}

The parallels $C_N^+$, $C_N^-$, $C_S^-$, $C_S^+$ separate the torus into four regions (see Fig.~\ref{Fig8}, left),  where the metric has constant signature:
\begin{equation}
\M_N \colon \Bigl|y-\frac{\pi}{2}\Bigr|< \frac{\pi}{4}, \ \ \,
\M_S \colon \Bigl|y-\frac{3\pi}{2}\Bigr|< \frac{\pi}{4},  \ \ \,
\M_E^+ \colon |y|< \frac{\pi}{4},  \ \ \,
\M_E^- \colon |y-\pi|< \frac{\pi}{4}.
\label{28}
\end{equation}
The North and the South regions $\M_N, \M_S$ are Riemannian,
while the outer and inner equatorial regions $\M_E^+, \M_E^-$ are Lorenzian.

Geodesics on the torus can be analyzed similarly to what was done for sphere (the main novelty is that the torus has regions with negative Gaussian curvature). This gives the following results:

\begin{itemize}
\item
The outer (large) equator $E^+$ and the inner (small) equator $E^-$ are only parallels of the torus being geodesics. Here the constant $\ha^2 = (\rho \pm 1)^2$ for $E^{\pm}$ respectively.

\item
$E^+$ is a unique geodesics in $\M^+$ that does not approach the set of parabolic points arbitrarily closely.
On the contrary, in $\M^-$ there exist an infinite number of geodesics possessing this property, they have the form of wave-like curves oscillating around $E^-$ (Fig.~\ref{Fig10}, left).

\item
All regions \eqref{28} contain bows of meridians, which are
timelike geodesics in $\M_N, \M_S$ and spacelike geodesics in $\M_E^+, \M_E^-$ ($\ha^2=0$).
\end{itemize}

Further we describe remaining geodesics (i.e., not being the equators $E^{\pm}$ or meridians) in three regions:
$\M_E^+$, $\M_E^-$, $\M_N$. Geodesics in the region $\M_S$ can be obtained from geodesics in $\M_N$ due to mirror symmetry with respect to the equatorial plane.

\begin{figure}[h]
\begin{center}
\includegraphics{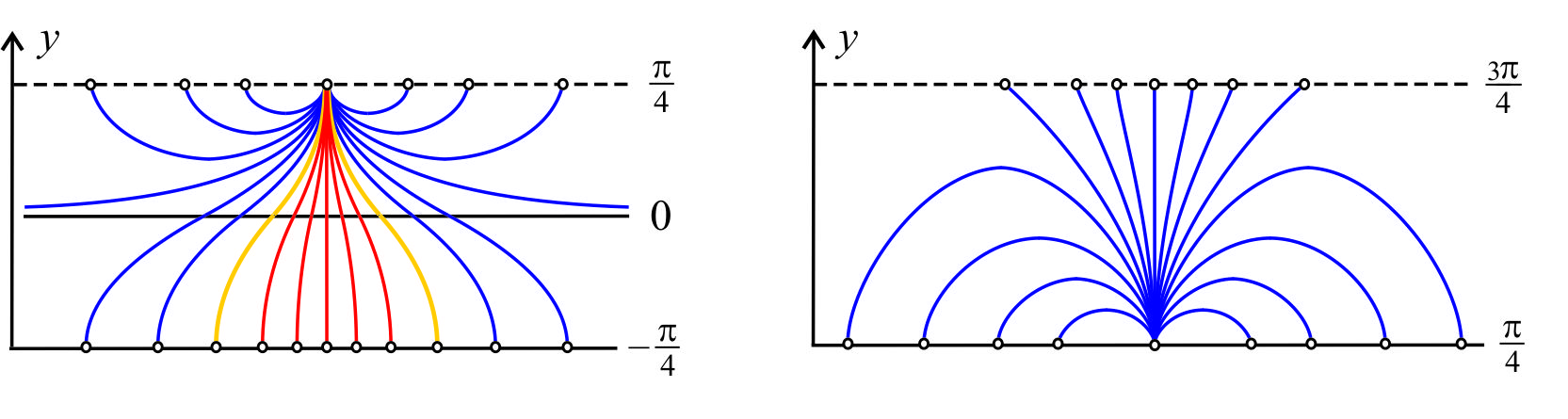}
\caption{ Geodesics outgoing from a parabolic point $q_0 \in C_N^+$
in the Lorenzian region $\M_E^+$ (left) and in the Riemannian region $\M_N$ (right).
}
\label{Fig9}
\end{center}
\end{figure}

\begin{figure}[h]
\begin{center}
\includegraphics{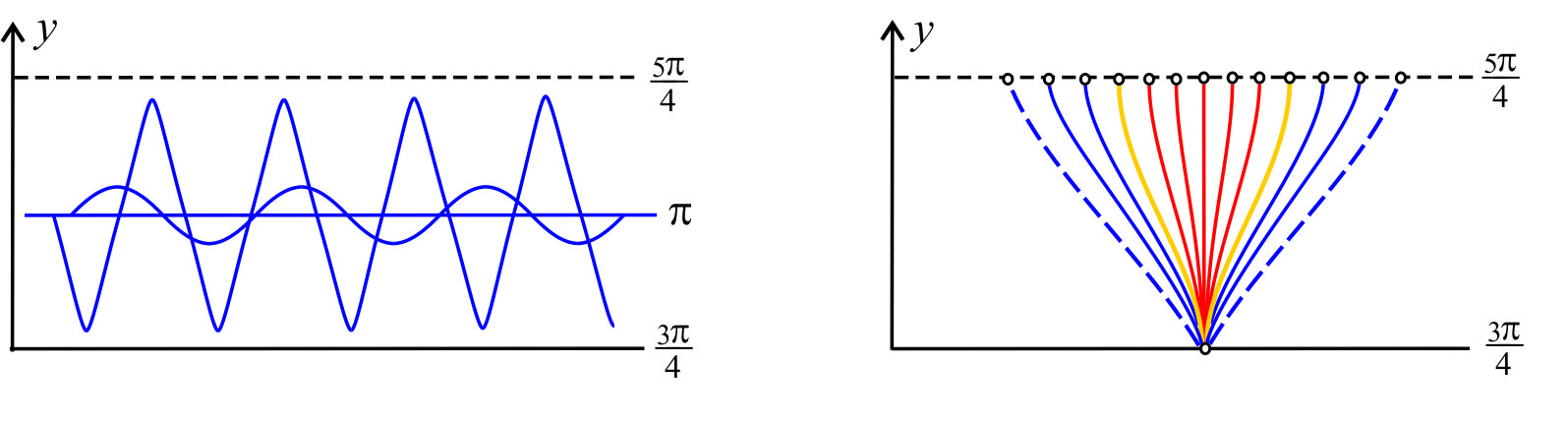}
\caption{
Geodesics in the Lorenzian region $\M_E^-$.
On the left: geodesics oscillating around the inner equator (class~1).
On the right: geodesics outgoing from a parabolic point $q_0 \in C_N^-$ (classes 2\,--\,5).
Timelike, spacelike, and isotropic geodesics are depicted as the blue, red, and yellow lines respectively.
Geodesics passing through $q_0$ with non-isotropic admissible directions (class~2) are depicted as the dashed lines.
}
\label{Fig10}
\end{center}
\end{figure}

\begin{itemize}
\item
Geodesics in the outer equatorial region $\M_E^+$ of the torus are similar to geodesics in the Lorenzian region of the sphere. They can be divided into seven classes similar to Table~2, with evident changes.\footnote{
The changes are as follows. In the second column the numbers 1 and 2 (except in $\ha^2$) should be
replaced by $(\rho+\frac{1}{\sqrt{2}})^2$ and $(\rho+1)^2$, respectively.
In all others columns the parallels $C_N$ and $C_S$ should be replaced by $C_N^+$ and $C_S^+$;
$N$- and $S$-hemispheres should be replaced by the North and South halves of the torus;
the equator $E$ should be replaced by the outer equator $E^+$.
}
The family of geodesics outgoing from a parabolic point $q_0 \in C_N^+$ in the region $\M_E^+$
is presented on Fig.~\ref{Fig9} (left).

\item
In the North region $\M_N$, all geodesics are curves of finite length with two endpoints.
The first endpoint is a cusp on $C_N^+$ with the isotropic tangential direction. 
The second endpoint can be a cusp on $C_N^+$ or a cusp on $C_N^-$ or a regular point on $C_N^-$.
In the first and second cases, the tangential direction is isotropic, while in the third case, it is one of two timelike admissible directions $p_{1,2}$; see Fig.~\ref{Fig9} (right).

\item
The inner equatorial region $\M_E^-$ contains an infinite number of timelike geodesics that oscillate around the inner equator $E^-$ and do not approach the boundary $C_N^- \cup C_S^-$ arbitrarily closely (Fig.~\ref{Fig10}, left). Also $\M_E^-$ contains an infinite number of geodesics joining the parallels $C_N^-$ and $C_S^-$, both endpoints of each such geodesic are cusps with the isotropic tangential directions or regular points with timelike admissible directions $p_{1,2}$ (Fig.~\ref{Fig10}, right).
Altogether, there are five classes enumerated in Table~3:
\end{itemize}

\begin{table}[h]
\label{Table3}
\begin{center}
\begin{tabular}{|c|c|c|c|c|c|}
\hline
{} & Type & Range of $\ha^2$ & 1st endpoint & 2nd endpoint &   Brief description  $\phantom{\Bigl|}$  \\
\hline
1  & timelike &  $(\rho - 1)^2 < \ha^2 < (\rho - \frac{1}{\sqrt{2}})^2$
$\phantom{\Bigl|}$ & --- & ---  & Oscillate around $E^-$ \\
\hline
2  & timelike &  $\ha^2 = (\rho - \frac{1}{\sqrt{2}})^2$
$\phantom{\Bigl|}$ & regular on $C_N^-$ & regular on $C_S^-$  &
Connect $C_N^-$ and $C_S^-$ \\
\hline
3  & timelike &  $(\rho - \frac{1}{\sqrt{2}})^2 < \ha^2 < +\infty$
$\phantom{\Bigl|}$ & cusp on $C_N^-$ & cusp on $C_S^-$  &
Connect $C_N^-$ and $C_S^-$ \\
\hline
4  & isotropic &  $\ha^2 = \infty$
$\phantom{\Bigl|}$ & cusp on $C_N^-$ & cusp on $C_S^-$  &
Connect $C_N^-$ and $C_S^-$ \\
\hline
5  & spacelike &  $-\infty < \ha^2 < 0$
$\phantom{\Bigl|}$ & cusp on $C_N^-$ & cusp on $C_S^-$  &
Connect $C_N^-$ and $C_S^-$ \\
\hline
\end{tabular}
\caption{
Five classes of geodesics (with the exception of the equator $E^-$ and meridians) in $\M_E^-$.
}
\end{center}
\end{table}

\subsection{Klein type and Grushin type metrics}\label{Klein-and-Grushin}

In this section, we consider geodesics in discontinuous metrics of the following two types:
\begin{equation}
ds^2= \frac{v(y)\,dx^2 +  w(y)\,dy^2}{y^2},  \ \quad   ds^2= \frac{v(y)}{y^2}\,dx^2 + w(y)\,dy^2 ,
\label{X1}
\end{equation}
where $v(y),w(y)$ are smooth positive functions.
The metrics \eqref{X1} are natural generalizations of the Klein metric \eqref{21} and the Grushin metric, respectively. Grushin type metrics have various applications in control theory and mechanics, see e.g.~\cite{BC}.

Both metrics \eqref{X1} have discontinuity on the line $\AA \colon y=0$ that lead to the same phenomenon as in the case of pseudo-Riemannian metrics: geodesics cannot pass through a point $q_0 \in \AA$ in arbitrary tangential directions, but only in the admissible direction $p=\infty$, see \cite{GR, Rem-Klein}.
Without loss of generality assume that $q_0=0$ (the origin) and $v(0)=w(0)=1$.
As proved in \cite{GR, Rem-Klein}, geodesics of the family $\Gamo$ have the form
\begin{equation}
x = \alf y^2 + o(y^2), \quad x = \alf y^3 + o(y^3), \quad \alf \in \bR,
\label{X11}
\end{equation}
for Klein type and Grushin type metrics, respectively.\footnote{
One of the differences among Klein type and Grushin type metrics can be found in the different asymptotic characters of the natural parametrization of geodesics tending to a point $q_0 \in \AA$. Naturally parametrized geodesics in Klein type metrics reach $q_0$ in infinite time, while in Grushin type metrics they reach it in finite time.
It follows from \eqref{3} similarly to what was done earlier for pseudo-Riemannian metrics.
}

\medskip

Unparametrized geodesics in the metrics \eqref{X1} are solutions of the Euler--Lagrange equation \eqref{6} with the Lagrangians $\LL = \sqrt{v+wp^2}/y$, $\LL = \sqrt{v+w(yp)^2}/y$, respectively. The corresponding the energy integrals \eqref{20} read
\begin{equation}
H(y,p) = \frac{v}{\sqrt{y^2(v + wp^2)}},  \ \quad
H(y,p) = \frac{v}{\sqrt{y^2(v + w(yp)^2)}}.
\label{X12}
\end{equation}
To establish the set $\HH$ of all possible constants $\ha^2$, one can substitute \eqref{X11} and the corresponding expressions $p^{-1} = 2\alf y + o(y)$, $p^{-1} = 3\alf y^2 + o(y^2)$ in the functions $H^2(y,p)$ from \eqref{X12} and pass to the limit $y \to 0$. Since $v(0)=w(0)=1$, for Klein type metrics this yields
\begin{equation}
\ha^2 = \lim_{y \to 0} H^2(y,p) =
\lim_{y \to 0} \frac{v^2(y)}{y^2v(y) + y^2w(y) (2\alf y + o(y))^{-2}} = 4\alf^2,
\label{X2}
\end{equation}
hence $\HH = [0, +\infty)$.
Analogous reasonings give the same result for Grushin type metrics.

\medskip

To define the behavior of geodesics $\gamapm \in \Gamo$, it is sufficient to consider the equation $v(y)/y^2 = \ha^2$ for all possible values $\ha^2 \ge 0$. This gives the following conclusions:

\begin{itemize}
\item
Geodesics $\gam^{\pm}_{\alf}$ and $\gam^{\pm}_{-\alf}$ with the same superscript return (or do not return)
to the initial line $y=0$ simultaneously.

\item
The family $\Gamo$ contains at least one couple $\gamapm$ that do not return to the initial line $y=0$:
$\gam_0^{\pm}$ with $\ha^2=0$ are the halves of the $y$-axis.

\item
The family $\Gamo$ contains an infinite number of geodesics that return to the initial line $y=0$.
Moreover, for any $\eps>0$ there exists $A>0$ such that all $\gamapm \in \Gamo$ with $|\alf|>A$
do not leave the strip $|y| < \eps$ and return to $y=0$.

\item
Geodesics $\gamap \in \Gamo$ (respectively, $\gamam \in \Gamo$) with all $\alf \neq 0$ return to the line $y=0$
if and only if $v(y)/y^2 \to 0$ as $y \to +\infty$ (respectively, $y \to -\infty$) and the condition $yv'(y) - 2v(y) \neq 0$ holds for all $y>0$ (respectively, $y<0$).
\end{itemize}

{\bf Example 3.3.}
For $v(y) = 1$, we have $v(y)/y^2 \to 0$ as $y \to \pm \infty$ and $yv'(y) - 2v(y) \neq 0$ for all $y$.
Hence $\gam_0^{\pm} \in \Gamo$ are only two geodesics that do not return to the line $y=0$,
while all $\gamapm \in \Gamo$ with $\alf \neq 0$ return.
In the case of the Klein metric \eqref{21} this statement is obvious: $\gamapm \in \Gamo$ with $\alf \neq 0$ are the halves of the circles $(x-x_*)^2 + y^2 = h_{\alf}^{-2}$.
\hfill\ebox

\begin{figure}[h!]
\begin{center}
\includegraphics{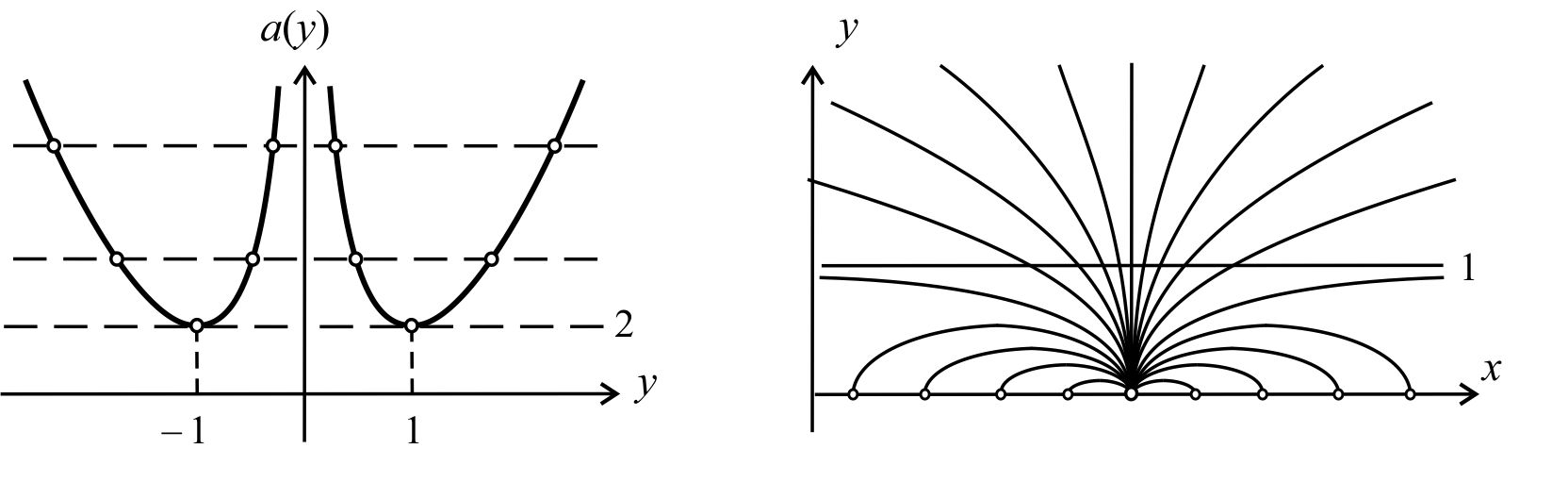}
\caption{
Illustration for Example~3.4.
On the left: the graph of the function $a(y) = v(y)/y^2 = y^2 + y^{-2}$.
On the right: geodesics $\gamap \in \Gamo$ for a point $q_0 \in \AA$.
}
\label{Fig11}
\end{center}
\end{figure}

\begin{figure}[h!]
\begin{center}
\includegraphics{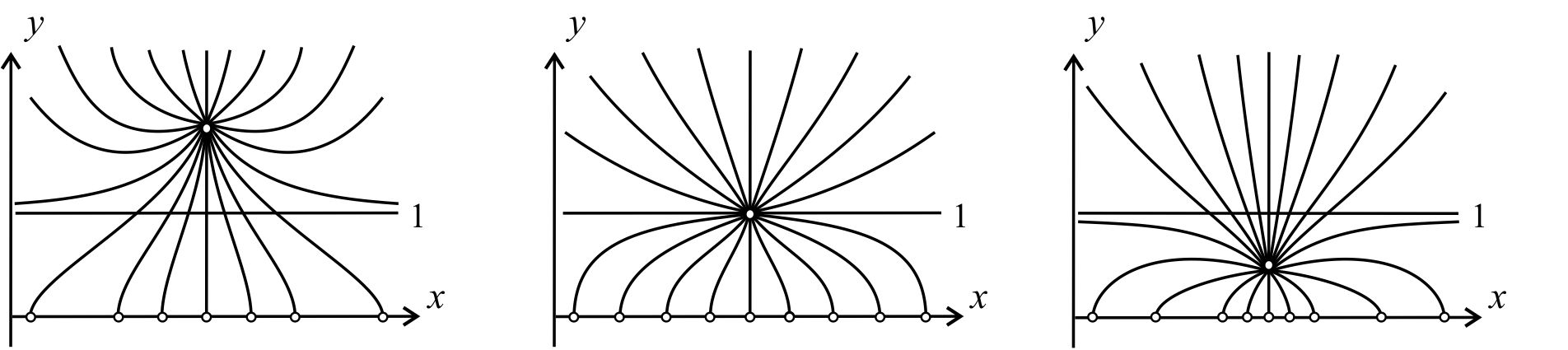}
\caption{
Illustration for Example~3.4.
Geodesics $\gamapm \in \Gamo$ for a point $q_0$ with $y_0>1$, $y_0=1$, $0<y_0<1$ (from left to right).
}
\label{Fig12}
\end{center}
\end{figure}

{\bf Example 3.4.}
Let $v(y) = 1 + y^4$.
The graph of the function $a(y) = v(y)/y^2$ is presented on Fig.~\ref{Fig11} (left).
It has two global minimums at $y = \pm 1$, hence there exist two horizontal geodesics $y = \pm 1$.
Further without loss of generality consider geodesics in the upper semiplane $y>0$.
Moreover, it is sufficient to analyze the family $\Gamo$ for $q_0 \in \AA$ ($y_0=0$) and for three arbitrary non-sungular points $q_0$ with $y_0>1$, $y_0=1$, $0<y_0<1$. In the last three cases, one can use results obtained in Section~\ref{Riemannian} for Riemannian metrics. The families $\Gamo$ for $q_0$ with $y_0>1$, $y_0=1$, $0<y_0<1$ are presented on Fig.~\ref{Fig12}.

Consider the family $\Gamo$ for $q_0 \in \AA$. There are three different classes of geodesics $\gamap \in \Gamo$ (and the same for $\gamam \in \Gamo$), presented on Fig.~\ref{Fig11} (right).
Firstly, $\gamap \in \Gamo$ with $\ha^2 > 2$, which do not leave the strip $0<y<1$ and return to the initial line $y=0$. Secondly, two geodesics $\gamap \in \Gamo$ with $\ha^2 = 2$ ($\alf = \pm \frac{1}{\sqrt{2}}$),
which do not leave the strip $0<y<1$, do not return to $y=0$, and asymptotically tend to the horizontal geodesic $y=1$. 
Finally, $\gamap \in \Gamo$ with $0 \le \ha^2 < 2$, which do not return to $y=0$ (the $y$-coordinate along them monotonically increases from $0$ to $+\infty$).
\hfill\ebox

\end{document}